# VAST VOLATILITY MATRIX ESTIMATION FOR HIGH-FREQUENCY FINANCIAL DATA

By Yazhen Wang[1] and Jian Zou

*University of Wisconsin-Madison and University of Connecticut*

High-frequency data observed on the prices of financial assets are commonly modeled by diffusion processes with micro-structure noise, and realized volatility-based methods are often used to estimate integrated volatility. For problems involving a large number of assets, the estimation objects we face are volatility matrices of large size. The existing volatility estimators work well for a small number of assets but perform poorly when the number of assets is very large. In fact, they are inconsistent when both the number, $p$, of the assets and the average sample size, $n$, of the price data on the $p$ assets go to infinity. This paper proposes a new type of estimators for the integrated volatility matrix and establishes asymptotic theory for the proposed estimators in the framework that allows both $n$ and $p$ to approach to infinity. The theory shows that the proposed estimators achieve high convergence rates under a sparsity assumption on the integrated volatility matrix. The numerical studies demonstrate that the proposed estimators perform well for large $p$ and complex price and volatility models. The proposed method is applied to real high-frequency financial data.

**1. Introduction.** Intra-day data observed on the prices of financial assets are often referred to as high-frequency financial data. Advances in technology make high-frequency financial data widely available nowadays for a host of different financial instruments on markets of all locations and at various scales, from individual bids to buy and sell, to the full distribution of such bids. The wide availability, in turn, stirs up an increasing demand for better modeling and statistical inference regarding the price and

Received October 2008; revised June 2009.

[1]Supported in part by NSF Grant DMS-05-04323, and this material was based on work supported by the National Science Foundation while author was working at the Foundation as a Program Director.

volatility dynamics of the assets. Diffusion processes are often employed to model high-frequency financial data, and various methodologies have been developed in past several years to estimate integrated volatility (or diffusion variance) over a period of time, say, a day. For a single asset, estimators of integrated volatility include realized volatility (RV) [Andersen et al. (2003), Barndorff-Nielsen and Shephard (2002)], bi-power realized variation (BPRV) [Barndorff-Nielsen and Shephard (2006)], two-time scale realized volatility (TSRV) [Zhang, Mykland and Aït-Sahalia (2005)], multiple-time scale realized volatility (MSRV) [Zhang (2006)], wavelet realized volatility (WRV) [Fan and Wang (2007)], kernel realized volatility (KRV) [Barndorff-Nielsen et al. (2008a)], pre-averaging realized volatility [Jacod et al. (2007)] and Fourier realized volatility (FRV) [Mancino and Sanfelici (2008)]. For multiple assets, we encounter a so-called non-synchronization problem which refers to as the fact that transactions for different assets often occur at distinct times, and the high-frequency prices of the assets are recorded at mismatched time points. Hayashi and Yoshida (2005) and Zhang (2007) have proposed to estimate integrated co-volatility of two assets based on overlap intervals and previous ticks. Barndorff-Nielsen and Shephard (2004) considered estimation of integrated co-volatility for synchronized high-frequency data.

A large number of assets are usually involved with in asset pricing, portfolio allocation, and risk management. One key problem we face is to estimate an integrated volatility matrix of large size for the assets. The scenario fits in to the so-called small $n$ and large $p$ or large $n$ but much larger $p$ problem, a current hot topic in statistics. The existing volatility estimation methods work well only for the cases of a single asset or a small number of assets where volatility is either scalar or a small matrix. Their poor behaviors for a large volatility matrix are indicated by random matrix theory and large covariance matrix estimation. Although idealized, the following example is still able to illustrate the point. Consider $p$ assets over unit time interval with all prices following independent Black-Scholes model with zero drift and unit volatility. Then the log prices obey independent standard Brownian motions, and the true integrated volatility matrix $\boldsymbol{\Gamma}$ is equal to the identity matrix $\mathbf{I}_p$. Assume that we observe all $p$ prices without noise at the same time grids $t_\ell = \ell/n$ for $\ell = 0, 1, \ldots, n$. The corresponding returns are i.i.d. normal random variables with mean zero and variance $1/n$. For this case, the existing best estimator of $\boldsymbol{\Gamma}$ is the RV, $\widehat{\boldsymbol{\Gamma}}$, with the following representation:

$$\widehat{\boldsymbol{\Gamma}} = (\widehat{\Gamma}_{ij}), \qquad \widehat{\Gamma}_{ij} = \frac{1}{n} \sum_{\ell=1}^{n} Z_{i\ell} Z_{j\ell} \text{ for } 1 \leq i, j \leq p$$

where $Z_{i\ell}$, $\ell = 1, \ldots, n$, $i = 1, \ldots, p$, are independent $N(0,1)$ random variables. It is widely known that $\widehat{\boldsymbol{\Gamma}}$ is a poor estimator of $\boldsymbol{\Gamma}$ when both $n$ and $p$ are large [Johnstone ([2001](#)), Johnstone and Lu ([2009](#)), El Karoui ([2007](#), [2008](#)), Bickel and Levina ([2008a](#), [2008b](#))]. In fact, for $n$ and $p$ both going to infinity but $p/n \to c$, the largest eigenvalue of $\widehat{\Gamma}$ asymptotically behaves like $(1 + \sqrt{c})^2$ while all true eigenvalues are equal to 1.

This paper develops a methodology for estimating large volatility matrices based on high-frequency financial data. We establish asymptotic theory for the proposed estimators under sparsity or decay assumptions on integrated volatility matrices as both $n$ and $p$ go to infinity. The estimators proposed in this paper are constructed as follows. In stage one we select grids as pre-sampling frequencies, construct a realized volatility matrix using previous tick method according to each pre-sampling frequency and then take the average of the constructed realized volatility matrices as the stage one estimator, which we call the ARVM estimator. In stage two we regularize the ARVM estimator to yield good consistent estimators of the large integrated volatility matrix. We consider two regularizations: thresholding and banding, developed by Bickel and Levina ([2008a](#), [2008b](#)) in the context of covariance matrix estimation. Thresholding a matrix is to retain only the elements whose absolute values exceed a given value which is called threshold and to replace the others by zero. Thresholding technique was introduced in wavelet literature for function estimation and image analysis [Wang ([2006](#))] where a function is known to have sparse representations in the sense that there are a relatively small number of important terms in its representations, but neither the number nor the locations of the important terms are known. We use thresholding to pick up the important terms for constructing estimators. For a sparse matrix, the small number of elements with large values are important. We need to locate those elements and estimate their values. The thresholding ARVM estimator is designed to find its elements of large magnitude along with their locations. "Banding" a matrix is to keep only the elements in a band along its diagonal and replace others by zero. Banding is analog to smoothing in nonparametric function estimation where in Taylor or orthogonal expansions of a given function, the locations of important terms are known. We simply choose the important terms for building estimators. For a matrix with elements decaying away from its diagonal, important terms are elements within a band along the diagonal, and banding ARVM estimator is used to select its elements within the band. The regularized ARVM estimators provide better volatility estimation that can greatly enhance portfolio allocation and risk management. With the volatility matrix estimators obtained from high-frequency data, we are able to investigate volatility dynamics directly and significantly improve volatility forecasting. See Andersen, Bollerslev and Diebold ([2004](#)) and Wang, Yao and Zou ([2008](#)).

We have shown that for a sparse integrated volatility matrix, the thresholded ARVM estimator not only consistently estimates the integrated volatility matrix but also achieves a high convergence rate when $p$ is allowed to grow as fast as a power of sample size with the power depending on the number of moments imposed on volatility processes and micro-structure noise. When it is known that the integrated volatility matrix has elements decaying away from its diagonal, the banded ARVM estimator is consistent and may enjoy a slightly higher convergence rate than the thresholded ARVM estimator. We have conducted extensive simulation studies for sophisticated price and volatility models with large $p$. The simulation studies demonstrate that the proposed estimators perform well for finite $p$ and sample size. The proposed method is applied to high-frequency data on 630 stocks traded in the Shanghai Stock Exchange.

The problem considered in this paper is much more complex than covariance matrix estimation, and our technical analyses rely on delicate treatments of diffusion processes and noise. Consequently the assumptions imposed and convergence rates obtained are different from those in the covariance matrix context. First, because of micro-structure noise in high-frequency financial data, true assets' prices are not directly observable; second, observations for continuous price processes are available only at discrete time points; third, price data on multiple assets are nonsynchronized; fourth, randomness in the observed data is caused by micro-structure noise as well as uncertainties in price and volatility processes. Because of noise, nonsynchronization and discretization, the convergence rates for volatility matrix estimation depend on sample size through a rate slower than the usual square root rate for covariance matrix estimation. Due to the multiple random sources in the data, for commonly used price and volatility models, the observed data cannot be sub-Gaussian. As a result, the convergence rates increase in $p$ faster for volatility matrix estimation than for covariance matrix estimation.

The rest of the paper is organized as follows. Section 2 introduces the basic models for high-frequency data and the estimation problem. The proposed methodology is presented in Section 3. The asymptotic theory is established in Section 4. Numerical studies are reported in Section 5. All proofs are relegated in Sections 6 and 7.

## 2. The set-up.

2.1. *Observed data.* Consider $p$ assets, and let $X_i(t)$ be the true log price at time $t$ of the $i$th asset, $i = 1, \ldots, p$. Suppose that we have high-frequency data for which the true log price of the $i$th asset is observed at times $t_{i\ell}$, $\ell = 1, \ldots, n_i$, and denote by $Y_i(t_{i\ell})$ the observed log price of the $i$th asset at time $t_{i\ell}$. Because of the nonsynchronization problem, typically $t_{i\ell} \neq t_{j\ell}$ for any $i \neq j$. The high-frequency data are usually contaminated with micro-structure noise in the sense that the observed log price $Y_i(t_{i\ell})$ is a noisy version of the corresponding true log price $X_i(t_{i\ell})$. It is common to assume

$$(1) \qquad Y_i(t_{i\ell}) = X_i(t_{i\ell}) + \varepsilon_i(t_{i\ell}), \qquad i = 1, \dots, p, \ell = 1, \dots, n_i,$$

where $\varepsilon_i(t_{i\ell})$, $i = 1, \dots, p$, $\ell = 1, \dots, n_i$ are independent noises with mean zero, for each fixed $i$, $\varepsilon_i(t_{i\ell})$, $\ell = 1, \dots, n_i$ are i.i.d. random variables with variance $\eta_i$, and $\varepsilon_i(\cdot)$ and $X_i(\cdot)$ are independent.

In the realized volatility literature, it is often assumed that micro-structure noise is i.i.d. and independent of the underlying price process. The simplistic assumption is used to study the effect of micro-structure noise on the volatility estimation. Recently Hansen and Lunde ([2006](#)) and Kalnina and Linton ([2008](#)), among others, have considered univariate micro-structure models where micro-structure noise has serial correlation and is correlated with the underlying price process.

2.2. *Price model.* Let $\mathbf{X}(t) = (X_1(t), \dots, X_p(t))^T$ be the vector of the true log prices at time $t$ of $p$ assets. Following finance theory we assume that $\mathbf{X}(t)$ obeys a continuous-time diffusion model,

$$(2) \qquad d\mathbf{X}(t) = \boldsymbol{\mu}_t \, dt + \boldsymbol{\sigma}_t^T \, d\mathbf{B}_t, \qquad t \in [0, 1],$$

where $\boldsymbol{\mu}_t = (\mu_1(t), \dots, \mu_p(t))^T$ is a drift vector, $\mathbf{B}_t = (B_{1t}, \dots, B_{pt})^T$ is a standard $p$-dimensional Brownian motion [i.e., $B_{it}$ are independent standard Brownian motions] and $\boldsymbol{\sigma}_t$ is a $p$ by $p$ matrix. The volatility of $\mathbf{X}(t)$ is given by

$$\boldsymbol{\gamma}(t) = (\gamma_{ij}(t))_{1 \leq i,j \leq p} = \boldsymbol{\sigma}_t^T \boldsymbol{\sigma}_t,$$

and its quadratic variation is equal to

$$[\mathbf{X}, \mathbf{X}]_t = \int_0^t \boldsymbol{\gamma}(s) \, ds = \left( \int_0^t \gamma_{ij}(s) \, ds \right)_{1 \leq i,j \leq p}, \qquad t \in [0, 1].$$

Our goal is to estimate the integrated volatility matrix,

$$\boldsymbol{\Gamma} = (\Gamma_{ij})_{1 \leq i,j \leq p} = \int_0^1 \boldsymbol{\gamma}(t) \, dt = \left( \int_0^1 \gamma_{ij}(t) \, dt \right)_{1 \leq i,j \leq p},$$

based on noisy and nonsynchronized observations $Y_i(t_{i\ell})$, $\ell = 1, \dots, n_i$, $i = 1, \dots, p$.

### 3. Estimation methodology.

3.1. *Realized volatility matrix.* Fix an integer $m$ and take $\tau_r, r = 1, \ldots, m$, to be the pre-determined sampling frequency. Let $\boldsymbol{\tau} = \{\tau_r, r = 1, \ldots, m\}$. For asset $i$, define previous-tick times

$$\tau_{i,r} = \max\{t_{i\ell} \le \tau_r, \ell = 1, \ldots, n_i\}, \qquad r = 1, \ldots, m.$$

Based on $\boldsymbol{\tau}$ we define realized co-volatility between assets $i$ and $j$ by

$$(3) \qquad \widehat{\Gamma}_{ij}(\boldsymbol{\tau}) = \sum_{r=1}^{m} [Y_i(\tau_{i,r}) - Y_i(\tau_{i,r-1})][Y_j(\tau_{j,r}) - Y_j(\tau_{j,r-1})],$$

and realized volatility matrix by

$$(4) \qquad \widehat{\boldsymbol{\Gamma}}(\boldsymbol{\tau}) = (\widehat{\Gamma}_{ij}(\boldsymbol{\tau})).$$

The pre-determined sampling frequency $\boldsymbol{\tau}$ is usually selected as regular grids. For a fixed $m$, there are $K = [n/m]$ classes of nonoverlap regular grids given by

$$\boldsymbol{\tau}^k = \{r/m, r = 1, \ldots, m\} + (k-1)/n = \{r/m + (k-1)/n, r = 1, \ldots, m\},$$
$$(5)$$

where $k = 1, \ldots, K$ and $n$ is the average sample size

$$n = \frac{1}{p} \sum_{i=1}^{p} n_i.$$

For each sampling frequency $\boldsymbol{\tau}^k$, using (3) and (4) we define realized co-volatility $\widehat{\Gamma}_{ij}(\boldsymbol{\tau}^k)$ between assets $i$ and $j$ and realized volatility matrix $\widehat{\boldsymbol{\Gamma}}(\boldsymbol{\tau}^k)$. Define

$$(6) \qquad \widehat{\Gamma}_{ij} = \frac{1}{K} \sum_{k=1}^{K} \widehat{\Gamma}_{ij}(\boldsymbol{\tau}^k), \qquad \widehat{\boldsymbol{\Gamma}} = (\widehat{\Gamma}_{ij}) = \frac{1}{K} \sum_{k=1}^{K} \widehat{\boldsymbol{\Gamma}}(\boldsymbol{\tau}^k).$$

Like TSRV in the univariate case [Zhang, Mykland and Aït-Sahalia (2005)], $\widehat{\boldsymbol{\Gamma}}$ is the average of $K$ realized volatility matrices $\widehat{\boldsymbol{\Gamma}}(\boldsymbol{\tau}^k)$. We use $\boldsymbol{\tau}^k$ to subsample data for computing $\widehat{\boldsymbol{\Gamma}}(\boldsymbol{\tau}^k)$ and then take their average to define $\widehat{\boldsymbol{\Gamma}}$. The purpose of subsampling and averaging is to handle noise and yield a better estimator.

We need to adjust the diagonal elements of $\widehat{\boldsymbol{\Gamma}}$ to account for the noise variances. Set $\boldsymbol{\eta} = \text{diag}(\eta_1, \ldots, \eta_p)$ where $\eta_i$ is the variance of noise $\varepsilon_i(t_{i\ell})$. We estimate $\eta_i$ by

$$(7) \qquad \widehat{\eta}_i = \frac{1}{2n_i} \sum_{\ell=1}^{n_i} [Y_i(t_{i,\ell}) - Y_i(t_{i,\ell-1})]^2,$$

and denote by $\widehat{\boldsymbol{\eta}} = \mathrm{diag}(\widehat{\eta}_1, \ldots, \widehat{\eta}_p)$ the estimator of $\boldsymbol{\eta}$. Define an estimator of $\boldsymbol{\Gamma}$ by

$$(8) \qquad \tilde{\boldsymbol{\Gamma}} = (\tilde{\Gamma}_{ij}) = \widehat{\boldsymbol{\Gamma}} - 2m\widehat{\boldsymbol{\eta}},$$

that is, we estimate element $\Gamma_{ij}$ of $\boldsymbol{\Gamma}$ by $\widehat{\Gamma}_{ij}$ for $i \neq j$ and $\widehat{\Gamma}_{ii} - 2m\widehat{\eta}_i$ for $i = j$. The diagonal elements of $\tilde{\boldsymbol{\Gamma}}$ are equal to TSRV of Zhang, Mykland and Aït-Sahalia (2005). We call $\tilde{\boldsymbol{\Gamma}}$ the averaging realized volatility matrix (ARVM) estimator.

High-frequency financial data are usually not equally spaced nor synchronized, and thus observations may be more or less dense for some assets than others or in some time intervals than others. An asset may have no observation between two consecutive time points in a sampling frequency; then the term involving these two consecutive time points in (3) is equal to zero, and from (6) we can see that the ARVM estimator automatically adjust for data with varying denseness.

3.2. *Regularize ARVM estimator.* For small $p$, $\tilde{\boldsymbol{\Gamma}}$ provides a good estimator for $\boldsymbol{\Gamma}$. However, $\tilde{\boldsymbol{\Gamma}}$ has a poor performance when $p$ gets very large. It is well known that even for constant $\boldsymbol{\gamma}(t)$, when $n$ and $p$ both go to infinity, estimators like $\tilde{\boldsymbol{\Gamma}}$ are inconsistent. In particular, when $p$ is very large, the eigenvalues and eigenvectors of $\tilde{\boldsymbol{\Gamma}}$ are far from those corresponding to $\boldsymbol{\Gamma}$ [see Bickel and Levina (2008a, 2008b), Johnstone (2001) and Johnstone and Lu (2009)].

We need to impose some structure on $\boldsymbol{\Gamma}$ and regularize $\tilde{\boldsymbol{\Gamma}}$ in order to estimate $\boldsymbol{\Gamma}$ consistently. Following Bickel and Levina (2008a, 2008b) we consider decay or sparsity assumptions on $\boldsymbol{\Gamma}$ and regularize $\tilde{\boldsymbol{\Gamma}}$ with banding or thresholding as follows.

Decay condition: We assume that the elements of $\boldsymbol{\Gamma}$ decay when moving away from its diagonal,

$$(9) \qquad |\Gamma_{ij}| \leq \frac{M}{1 + |i - j|^{\alpha+1}}, \qquad 1 \leq i, j \leq p, E[M] \leq C,$$

where $M$ is a positive random variable, and $C$ and $\alpha$ are positive generic constants.

Sparsity condition: We assume that $\boldsymbol{\Gamma}$ satisfies

$$(10) \qquad \sum_{j=1}^{p} |\Gamma_{ij}|^{\delta} \leq M\pi(p), \qquad i = 1, \ldots, p, E[M] \leq C,$$

where $M$ is a positive random variable, $0 \leq \delta < 1$, and $\pi(p)$ is a deterministic function of $p$ that grows very slowly in $p$.

Examples of $\pi(p)$ include 1, $\log p$ and a small power of $p$. The case of $\delta = 0$ in (10) corresponds so that each row of $\boldsymbol{\Gamma}$ has at most $M\pi(p)$ number of nonzero elements. Decay condition (9) corresponds to a special case of sparsity condition (10) with $\delta = 1/(\alpha + 1)$ and $\pi(p) = \log p$ or $1/(\alpha + 1) < \delta < 1$ and $\pi(p) = 1$.

The decay condition depends on the order of $p$ assets in the log price vector $\mathbf{X}(t)$. As stocks have no natural ordering, the decay condition may not hold for real volatility matrices of stock returns. As a result, for volatility matrix estimation in financial applications, sparsity is much more realistic than the decay assumption. Examples of sparse matrices include block diagonal matrices, matrices with decay elements from diagonal, matrices with relatively small number of nonzero elements in each row or column and matrices obtained by randomly permuting rows and columns of above matrices.

For $\mathbf{\Gamma}$ satisfying decay condition (9), its large elements are within a band along its diagonal, and the elements outside the band are negligible. We regularize $\tilde{\mathbf{\Gamma}}$ by banding, which is to keep only its elements in a band along its diagonal and replace others by zero. Specifically, the definition of banding $\tilde{\mathbf{\Gamma}}$ is given by

$$\mathcal{B}_b[\tilde{\mathbf{\Gamma}}] = (\tilde{\Gamma}_{ij} 1(|i - j| \le b)),$$

where $b$ is a banding parameter, and $1(|i - j| \le b)$ is the indicator of $\{(i, j), |i - j| \le b\}$. The $(i, j)$th element of $\mathcal{B}_b[\tilde{\mathbf{\Gamma}}]$ is equal to $\tilde{\Gamma}_{ij}$ for $|i - j| \le b$ and zero, otherwise. We call $\mathcal{B}_b[\tilde{\mathbf{\Gamma}}]$ the BARVM estimator.

If $\mathbf{\Gamma}$ satisfies sparsity condition (10), its important elements are those whose absolute values are above a certain threshold. We regularize $\tilde{\mathbf{\Gamma}}$ by thresholding which is to retain its elements whose absolute values exceed a given value and replace others by zero. Specifically, we define the thresholding of $\tilde{\mathbf{\Gamma}}$ by

$$\mathcal{T}_\varpi[\tilde{\mathbf{\Gamma}}] = (\tilde{\Gamma}_{ij} 1(|\tilde{\Gamma}_{ij}| \ge \varpi)),$$

where $\varpi$ is threshold. The $(i, j)$th element of $\mathcal{T}_\varpi[\tilde{\mathbf{\Gamma}}]$ is equal to $\tilde{\Gamma}_{ij}$ if its absolute value is greater or equal to $\varpi$ and zero, otherwise. We call $\mathcal{T}_\varpi[\tilde{\mathbf{\Gamma}}]$ TARVM estimator.

Like most of existing co-volatility matrix estimators, we cannot guarantee the positiveness of the ARVM estimator for finite sample. As the banding and thresholding procedures do not resolve the positiveness issue, the BARVM and TARVM estimators may not be positive for a finite sample. Recently Barndorff-Nielsen et al. (2008b) has developed a kernel-based method with refresh sampling time technique to produce a semi-positive co-volatility matrix estimator. The estimator is designed for fixed $p$ and must suffer from the same drawback as the ARVM estimator for large $p$; it will be interesting to apply the regularization procedures to the semi-positive matrix estimator and investigate their asymptotic behaviors for large $p$ and $n$.

Banding is analog to smoothing in nonparametric function estimation where in the representations of a target function by Taylor or orthogonal expansions the locations of important terms in the expansions are known. We simply select those important terms to keep for building estimators. For a matrix with decaying elements from its diagonal, important terms are elements within a band along the diagonal, and banding is used to select the elements within the band. Thresholding is utilized for estimating a function with sparse representations, where we know that there are a relatively small number of important terms in its representations, but neither the number nor the locations of the important terms are known. We use thresholding to pick up the important terms for constructing estimators. For a sparse matrix, all we know is that a relatively small number of the elements with large values essentially matter. We need to locate those elements and estimate their values. Thresholding is designed to find the elements of large magnitude and their locations.

For the BARVM and TARVM estimators, we need to select proper values for banding parameter $b$ and threshold $\varpi$ from data. Data-dependent methods for selecting $b$ and $\varpi$ are illustrated at the end of Section 5.3 for simulated data and at the beginning of Section 5.5 for real data.

**4. Asymptotic theory.** First we fix some notations for the theoretical analysis. Given a $p$-dimensional vector $\mathbf{x} = (x_1, \ldots, x_p)^T$ and a $p$ by $p$ matrix $\mathbf{U} = (U_{ij})$, define their $\ell_d$-norms as follows:

$$\|\mathbf{x}\|_d = \left(\sum_{i=1}^{p} |x_i|^d\right)^{1/d}, \qquad \|\mathbf{U}\|_d = \sup\{\|\mathbf{U}\mathbf{x}\|_d, \|\mathbf{x}\|_d = 1\}, \qquad d = 1, 2, \infty.$$

Then $\|\mathbf{U}\|_2$ is equal to the square root of the largest eigenvalue of $\mathbf{U}\mathbf{U}^T$,

$$\|\mathbf{U}\|_1 = \max_{1 \le j \le p} \sum_{i=1}^{p} |U_{ij}|, \qquad \|\mathbf{U}\|_\infty = \max_{1 \le i \le p} \sum_{j=1}^{p} |U_{ij}|,$$

and

$$\|\mathbf{U}\|_2^2 \le \|\mathbf{U}\|_1 \|\mathbf{U}\|_\infty.$$

For symmetric $\mathbf{U}$, $\|\mathbf{U}\|_2$ is equal to its largest absolute eigenvalue, and $\|\mathbf{U}\|_2 \le \|\mathbf{U}\|_1 = \|\mathbf{U}\|_\infty$.

Next we state some technical conditions.

A1: For some $\beta \ge 2$,

$$\max_{1 \le i \le p} \max_{0 \le t \le 1} E[|\gamma_{ii}(t)|^\beta] < \infty, \qquad \max_{1 \le i \le p} \max_{0 \le t \le 1} E[|\mu_i(t)|^{2\beta}] < \infty,$$

$$\max_{1 \le i \le p} E[|\varepsilon_i(t_{i\ell})|^{2\beta}] < \infty.$$

A2: Each of $p$ assets has at least one observation between $\tau_r^k$ and $\tau_{r+1}^k$. With $n = (n_1 + \cdots + n_p)/p$, we assume

$$C_1 \leq \min_{1 \leq i \leq p} \frac{n_i}{n} \leq \max_{1 \leq i \leq p} \frac{n_i}{n} \leq C_2, \qquad \max_{1 \leq i \leq p} \max_{1 \leq \ell \leq n_i} |t_{i\ell} - t_{i,\ell-1}| = O(n^{-1}),$$

$$m = o(n).$$

THEOREM 1. *Under models* (1) *and* (2) *and conditions* A1 *and* A2 *we have for all* $1 \leq i, j \leq p$,

$$(11) \qquad E(|\tilde{\Gamma}_{ij} - \Gamma_{ij}|^\beta) \leq C e_n^\beta,$$

*where* $C$ *is a generic constant free of* $n$ *and* $p$, *and the convergence rate* $e_n^\beta$ *given below is equal to the sum of terms with powers of* $n$ *and* $K = [n/m]$ *which depend on whether the observed data in the model specification have micro-structure noise or not.*

(1) *If there is micro-structure noise in model* (1),

$$e_n^\beta = (Kn^{-1/2})^{-\beta} + K^{-\beta/2} + (n/K)^{-\beta/2} + K^{-\beta} + n^{-\beta/2}.$$

*Thus with* $K \sim n^{2/3}$ *we have* $e_n \sim n^{-1/6}$.

(2) *If there is no micro-structure noise [i.e.* $\varepsilon_i(t_{i\ell}) = 0$ *and* $Y_i(t_{i\ell}) = X_i(t_{i\ell})$] *in model* (1),

$$e_n^\beta = (n/K)^{-\beta/2} + K^{-\beta} + n^{-\beta/2}.$$

*Thus with* $K \sim n^{1/3}$ *we have* $e_n \sim n^{-1/3}$.

REMARK 1. The convergence rate $e_n$ can be attributed to three sources due to noise, nonsynchronization and discrete observations for continuous process $\mathbf{X}(t)$. Because of micro-structure noise in high-frequency financial data, true log-price $\mathbf{X}(t)$ is not directly observable. Furthermore, as a continuous process, $\mathbf{X}(t)$ is observed with noise only at discrete time points. Consequently the convergence rate $e_n$ is slower than $n^{-1/2}$. In fact, the optimal convergence rate for the univariate noise case is $n^{-1/4}$; the nonsynchronization for multiple assets further complicates the problem. The terms in convergence rates $e_n^\beta$ given by Theorem 1 can be identified to associate with specific sources as follows. The terms $(Kn^{-1/2})^{-\beta} + K^{-\beta/2}$ in $e_n^\beta$ are due to noise with $K^{-\beta} + n^{-\beta/2}$ contributed by nonsynchronization. Because $\mathbf{X}(t)$ is observed at discrete time points, we need to discretize $\mathbf{X}(t)$ and use its discretization to approximate integrated volatility matrix. Term $(n/K)^{-\beta/2}$ in $e_n^\beta$ is attributed to the approximation error due to the discretization of $\mathbf{X}(t)$. These are clearly spelled out by Propositions 1–3 in Section 7 for the proof of Theorem 1. The contributions of the three sources to the convergence rates for TSRV, MSRV, realized co-volatility estimators have been shown in Zhang, Mykland and Aït-Sahalia (2005) and Zhang (2006, 2007).

THEOREM 2. *Assume that $\mathbf{\Gamma}$ satisfies sparsity condition (10). Then under models (1) and (2) and conditions* A1 *and* A2, *we have*

$$\|\mathcal{T}_\varpi[\tilde{\mathbf{\Gamma}}] - \mathbf{\Gamma}\|_2 \leq \|\mathcal{T}_\varpi[\tilde{\mathbf{\Gamma}}] - \mathbf{\Gamma}\|_\infty = O_P(\pi(p)[e_n p^{2/\beta} h_{n,p}]^{1-\delta}),$$

*where $e_n$ is given in Theorem 1, $\varpi = e_n p^{2/\beta} h_{n,p}$, and $h_{n,p}$ is any sequence converging to infinity arbitrarily slow with one example $h_{n,p} = \log \log(n \wedge p)$.*

REMARK 2. The convergence rate in Theorem 2 is nearly equal to $\pi(p)[e_n \times p^{2/\beta}]^{1-\delta}$. Since $e_n \sim n^{-1/6}$ for the noise case and $e_n \sim n^{-1/3}$ for the noiseless case, in order to make $e_n p^{2/\beta}$ go to zero, $p$ needs to grow more slowly than $n^{\beta/12}$ for the noise case and $n^{\beta/6}$ for the noiseless case.

THEOREM 3. *Assume that $\mathbf{\Gamma}$ satisfies decay condition (9). Then under models (1) and (2) and conditions* A1 *and* A2, *we have that*

$$\|\mathcal{B}_b[\tilde{\mathbf{\Gamma}}] - \mathbf{\Gamma}\|_2 \leq \|\mathcal{B}_b[\tilde{\mathbf{\Gamma}}] - \mathbf{\Gamma}\|_\infty = O_P([e_n p^{1/\beta}]^{\alpha/(\alpha+1+1/\beta)}),$$

*where we select banding parameter $b$ of order $(e_n p^{1/\beta})^{-1/(\alpha+1+1/\beta)}$.*

REMARK 3. For $\mathbf{\Gamma}$ satisfying decay condition (9), the sparsity condition is held with $\delta = 1/(\alpha + 1)$ and $\pi(p) = \log p$. The convergence rate corresponding to Theorem 2 under the sparsity condition has a leading factor of order $[e_n p^{2/\beta}]^{\alpha/(\alpha+1)}$. Comparing it with the rate in Theorem 3, we conclude that the two convergence rates are quite close for reasonably large $\beta$.

REMARK 4. The convergence rates in Bickel and Levina (2008a, 2008b) for large covariance matrix estimation increase in matrix size $p$ through a power of $\log p$, but the convergence rates in Theorems 2 and 3 grow with $p$ through a power of $p$. The slower convergence rates here are due to the intrinsic complexity of our problem. The $\log p$ factor in the convergence rates of covariance matrix estimation is attributed to Gaussianity or sub-Gaussianity imposed on the observed data. In our set-up, observations $Y_i(t_{i\ell})$ from model (1) have random sources from both micro-structure noise $\varepsilon_i(t_{i\ell})$ and true log price $\mathbf{X}(t)$ given by model (2). The term $\int_0^t \boldsymbol{\sigma}_s^T d\mathbf{B}_s$ in $\mathbf{X}(t)$ does not obey sub-Gaussianity for common price and volatility models. Even though we assume normality on $\varepsilon_i(t_{i\ell})$, the observed data $Y_i(t_{i\ell})$ cannot be sub-Gaussian for the price and volatility models. Consequently we employ realistic moment conditions in assumption A1, obtain convergence rates for the elements of $\tilde{\mathbf{\Gamma}}$ in Theorem 1 and derive subsequent convergence rates with a power of $p$ for the regularized $\tilde{\mathbf{\Gamma}}$ in Theorems 2 and 3.

REMARK 5. For Gaussian observations, Cai, Zhang and Zhou (2008) have established optimal convergence rates for estimating a covariance matrix which is assumed to belong to a class of matrices satisfying the decay condition. The convergence rate for the minimax risk based on the squared $\ell_2$ norm is equal to the minimum of $n^{-2\alpha/(2\alpha+1)} + \log p/n$ and $p/n$. The result indicates that the convergence rate in Bickel and Levina (2008a) is suboptimal. It is very interesting and challenging to find optimal convergence rates for the volatility matrix estimation problem in our set-up.

## 5. Numerical studies.

5.1. *High-frequency data.* The real data set for our numerical studies is high-frequency tick by tick price data on 630 stocks traded in the Shanghai Stock Exchange over 177 days in 2003. For each day, we computed the ARVM estimator corresponding to $\tilde{\Gamma}$ defined in Section 3 where the predetermined sampling frequencies were selected to correspond with 5 minute returns. This yielded 177 matrices of size 630 by 630 as ARVM estimators of integrated volatility matrices over the 177 days. The average of these 177 matrices was then evaluated and denoted by $\Theta$.

5.2. *The simulation model.* In our simulation study the true log price $\mathbf{X}(t)$ of $p$ assets is generated from model (2) with zero drift, namely, the diffusion model,

$$(12) \qquad d\mathbf{X}(t) = \boldsymbol{\sigma}_t^T \, d\mathbf{B}_t, \qquad t \in [0,1],$$

where $\mathbf{B}_t = (B_{1t}, \ldots, B_{pt})^T$ is a standard $p$-dimensional Brownian motion, and we take $\boldsymbol{\sigma}_t$ as a Cholesky decomposition of $\boldsymbol{\gamma}(t) = (\gamma_{ij}(t))_{1 \le i,j \le p}$ which is defined below. Given the diagonal elements of $\boldsymbol{\gamma}(t)$, we define its off-diagonal elements by

$$(13) \qquad \gamma_{ij}(t) = \{\kappa(t)\}^{|i-j|} \sqrt{\gamma_{ii}(t)\gamma_{jj}(t)}, \qquad 1 \le i \ne j \le p,$$

where process $\kappa(t)$ is given by

$$(14) \qquad \kappa(t) = \frac{e^{2u(t)} - 1}{e^{2u(t)} + 1}, \qquad du(t) = 0.03[0.64 - u(t)]\,dt + 0.118u(t)\,dW_{\kappa,t},$$

$$(15) \qquad W_{\kappa,t} = \sqrt{0.96}W_{\kappa,t}^0 - 0.2\sum_{i=1}^{p} B_{it}/\sqrt{p};$$

$W_{\kappa,t}^0$ is a standard 1-dimensional Brownian motion independent of $\mathbf{B}_t$. Model (14) is taken from Barndorff-Nielsen and Shephard (2002, 2004).

The diagonal elements of $\boldsymbol{\gamma}(t)$ are generated from four common stochastic volatility models with leverage effect. The four volatility processes are geometric Ornstein–Uhlenbeck processes, the sum of two CIR processes [Cox, Ingersoll and Ross (1985) and Barndorff-Nielsen and Shephard (2002)], the volatility process in Nelson's GARCH diffusion limit model [Wang (2002)] and two-factor log-linear stochastic volatility process [Huang and Tauhen (2005)]. Specifically, let $\mathbf{U}_t^1 = (U_{1t}^1, \ldots, U_{pt}^1)^T$ and $\mathbf{U}_t^2 = (U_{1t}^2, \ldots, U_{pt}^2)^T$ be two independent standard $p$-dimensional Brownian motions which are independent of $\mathbf{B}_t$ and $W_{\kappa,t}^0$, and then define two $p$-dimensional Brownian motions, $\mathbf{W}_t^1 = (W_{1t}^1, \ldots, W_{pt}^1)^T$ and $\mathbf{W}_t^2 = (W_{1t}^2, \ldots, W_{pt}^2)^T$, by

$$W_{it}^1 = \rho_i B_{it} + \sqrt{1 - \rho_i^2} U_{it}^1, \qquad W_{it}^2 = \rho_i B_{it} + \sqrt{1 - \rho_i^2} U_{it}^2, \tag{16}$$

where we choose the following negative values for $\rho_i$ to reflect the leverage effect,

$$\rho_i = \begin{cases} -0.62, & 1 \le i \le p/4, \\ -0.50, & p/4 < i \le p/2, \\ -0.25, & p/2 < i \le 3p/4, \\ -0.30, & 3p/4 < i \le p. \end{cases}$$

We generate $\gamma_{ii}(t)$ as follows.

(1) For $1 \le i \le p/4$, $\gamma_{ii}(t)$ are drawn from the geometric Ornstein–Uhlenbeck model driving by $W_{it}^1$ [Barndorff-Nielsen and Shephard (2002)],

$$d \log \gamma_{ii}(t) = -0.6(0.157 + \log \gamma_{ii}(t)) \, dt + 0.25 \, dW_{it}^1. \tag{17}$$

(2) For $p/4 < i \le p/2$, $\gamma_{ii}(t)$ are drawn from the sum of two CIR processes [Barndorff-Nielsen and Shephard (2002)],

$$\gamma_{ii}(t) = 0.98(v_{1,t} + v_{2,t}), \tag{18}$$

where $v_{1,t}$ and $v_{2,t}$ obey two CIR models driving by $W_{it}^1$ and $W_{it}^2$, respectively,

$$dv_{1,t} = 0.0429(0.108 - v_{1,t}) \, dt + 0.1539\sqrt{v_{1,t}} \, dW_{i,t}^1, \tag{19}$$

$$dv_{2,t} = 3.74(0.401 - v_{2,t}) \, dt + 1.4369\sqrt{v_{2,t}} \, dW_{i,t}^2. \tag{20}$$

(3) For $p/2 < i \le 3p/4$, $\gamma_{ii}(t)$ are drawn from the volatility process in Nelson's GARCH diffusion limit model driving by $W_{it}^1$ [Barndorff-Nielsen and Shephard (2002)],

$$d\gamma_{ii}(t) = [0.1 - \gamma_{ii}(t)] \, dt + 0.2\gamma_{ii}(t) \, dW_{it}^1. \tag{21}$$

(4) For $3p/4 < i \le p$, $\gamma_{ii}(t)$ are drawn from the two-factor log-linear stochastic volatility model driving by $W_{it}^1$ and $W_{it}^2$ [Huang and Tauhen (2005)],

$$\gamma_{ii}(t) = e^{-6.8753}\text{s-exp}(0.04v_{1,t} + 1.5v_{2,t} - 1.2), \tag{22}$$

where

$$dv_{1,t} = -0.00137v_{1,t} \, dt + dW_{i,t}^1,$$

$$dv_{2,t} = -1.386v_{2,t} \, dt + (1 + 0.25v_{2,t}) \, dW_{i,t}^2, \tag{23}$$

$$\text{s-exp}(u) = \begin{cases} e^u, & \text{if } u \le \log(8.5), \\ 8.5\{1 - \log(8.5) + u^2/\log(8.5)\}^{1/2}, & \text{if } u > \log(8.5). \end{cases}$$

With $\gamma_{ii}(t)$ generated from above stochastic differential equations, we multiply $\gamma_{ii}(t)$ by $1000\theta_i$ where $\theta_i$ are the ordered (from the largest to the smallest) diagonal elements of $\boldsymbol{\Theta}$ defined in Section 5.1 as the average of 177 daily ARVM estimators for the high-frequency data from the Shanghai market. The adjustment is to roughly match simulated $\gamma_{ii}(t)$ with the magnitudes of the diagonal elements of the ARVM estimators for the stock data.

Finally the high-frequency data $Y_i(t_{i\ell})$ are simulated from model (1) with noise $\varepsilon_i(t_{i\ell})$ drawing from independent normal distributions with mean zero and standard deviation of three choices: $0.002\sqrt{\theta_i}$, $0.127\sqrt{\theta_i}$ and $0.2\sqrt{\theta_i}$ which correspond to low, medium and high noise levels. The standard deviation is chosen to reflect the empirical fact that relative noise level found in high frequency data typically ranges from $0.001\%$ to $0.01\%$ with $0.001\%$ for individual stock and $0.01\%$ for stock index. In our simulated example, the average volatility is around $1000\theta_i$, and thus the three noise standard deviations are translated into $0.002\%$, $0.004\%$ and $0.065\%$ of the average volatility or relative noise level, respectively.

5.3. *The simulation procedure.* We need to simulate $n$ values for the price and volatility processes at $t_\ell = \ell/n$, $\ell = 1, \ldots, n$. The procedure begins with the generation of matrices $\boldsymbol{\gamma}(t_\ell)$. First we use normalized partial sums of i.i.d. standard normal random variables to simulate four independent Brownian motions, a standard one-dimensional Brownian motion $W^0_{\kappa,t_\ell}$ and three standard $p$-dimensional Brownian motions, $\mathbf{B}_{t_\ell}$, $\mathbf{U}^1_{t_\ell}$ and $\mathbf{U}^2_{t_\ell}$, and compute $W_{\kappa,t_\ell}$, $\mathbf{W}^1_{t_\ell}$ and $\mathbf{W}^2_{t_\ell}$ according to (15) and (16). We then use the Euler scheme to simulate $\kappa(t_\ell)$ from (14) with $W_{\kappa,t_\ell}$ and $\gamma_{ii}(t_\ell)$ from (17)–(23) with corresponding components of $\mathbf{W}^1_{t_\ell}$ and $\mathbf{W}^2_{t_\ell}$. With available $\kappa(t_\ell)$ and $\gamma_{ii}(t_\ell)$, from (13) we evaluate off-diagonal elements $\gamma_{ij}(t_\ell)$, $i \neq j$. To speed up the simulation of $\boldsymbol{\gamma}(t_\ell) = (\gamma_{ij}(t_\ell))$, we have utilized the following tricks in R programming: (i) group all $p$ diagonal elements $\gamma_{11}(t_\ell), \ldots, \gamma_{pp}(t_\ell)$ into four vectors of dimension $p/4$ with each vector drawing from the same volatility model and update each vector in the Euler scheme; (ii) calculate the matrix product of column vector $(\gamma_{11}(t_\ell), \ldots, \gamma_{pp}(t_\ell))^T$ and row vector $(\gamma_{11}(t_\ell), \ldots, \gamma_{pp}(t_\ell))$ and then take the element by element square root of the obtained matrix; (iii) utilize the Toeplitz matrix operation to evaluate $\kappa(t_\ell)^{i-j}$; (iv) use the matrix operation to compute element by element multiplication of the two matrices yielded in steps (ii) and (iii).

We take $\boldsymbol{\sigma}_{t_\ell}$ as a Cholesky decomposition of $\boldsymbol{\gamma}(t_\ell)$ and compute true logprice $\mathbf{X}(t_\ell)$ by

$$\mathbf{X}_{t_\ell} = \mathbf{X}_{t_{\ell-1}} + [\boldsymbol{\sigma}_{t_{\ell-1}}]^T[\mathbf{B}_{t_\ell} - \mathbf{B}_{t_{\ell-1}}].$$

Finally, data $Y_i(t_\ell)$, $i = 1, \ldots, p$, $\ell = 1, \ldots, n$, are obtained by adding to $X_i(t_\ell)$ normal noise $\epsilon_i(t_\ell)$ of mean zero and standard deviation $0.002\sqrt{\theta_i}$, $0.127\sqrt{\theta_i}$ and $0.2\sqrt{\theta_i}$ for the cases of low, medium and high noise levels, respectively.

After the simulation has generated volatility matrices $\boldsymbol{\gamma}(t_\ell)$, log-price values $X_i(t_\ell)$ and synchronized noisy data $Y_i(t_\ell)$, $i = 1, \ldots, p$, $t_\ell = \ell/n$, $\ell = 1, \ldots, n$, we numerically evaluate integrated volatility matrix $\boldsymbol{\Gamma}$ by $\sum_{\ell=1}^{n} \boldsymbol{\gamma}(t_\ell)/n$, compute ARVM estimator $\tilde{\boldsymbol{\Gamma}}$ according to (8) and calculate its banding $\mathcal{B}_b[\tilde{\boldsymbol{\Gamma}}]$ and thresholding $\mathcal{T}_\varpi[\tilde{\boldsymbol{\Gamma}}]$ as described in Section 3.2. Note that there is no need to store matrices $\boldsymbol{\gamma}(t_\ell)$ in the programming loop of $\ell = 1, \ldots, n$, because at each $\ell$ step of the loop all we need is to record $X_i(t_\ell)$ and $Y_i(t_\ell)$ and update the partial sum $\sum_{r=1}^{\ell} \boldsymbol{\gamma}(t_r)$ to the current step for the purpose of evaluation of $\boldsymbol{\Gamma}$ by the average of $\boldsymbol{\gamma}(t_\ell)$ at the end of the loop. Doing so will save huge computer storage space and prevent the simulation program from running out of computer memory.

We repeat the whole simulation procedure 500 times. The mean square error (MSE) of a matrix estimator is computed by averaging $\ell_2$-norms of the differences between the estimator and $\boldsymbol{\Gamma}$ over 500 repetitions. We use the MSEs of $\tilde{\boldsymbol{\Gamma}}$, $\mathcal{B}_b[\tilde{\boldsymbol{\Gamma}}]$ and $\mathcal{T}_\varpi[\tilde{\boldsymbol{\Gamma}}]$ to evaluate their performances. For estimators $\mathcal{B}_b[\tilde{\boldsymbol{\Gamma}}]$ and $\mathcal{T}_\varpi[\tilde{\boldsymbol{\Gamma}}]$, we select the values of $b$ and $\varpi$ by minimizing their respective MSEs.

We generate nonsynchronized data as follows. Instead of generating observations for the processes at $n$ time points, we simulate $\boldsymbol{\gamma}(t_\ell)$, $X_i(t_\ell)$, $Y_i(t_\ell)$, $i = 1, \ldots, p$, at $3n$ time points $t_\ell = \ell/(3n)$, $\ell = 1, \ldots, 3n$. Grouping together three consecutive time points we divide the $3n$ time points $t_\ell$ into $n$ groups $\{t_{3r-2}, t_{3r-1}, t_{3r}\}$, $r = 1, \ldots, n$. For the $i$th asset, we select one time point at random from each group; from the simulated $3n$ values of $Y_i(t_\ell)$ we choose $n$ values corresponding to the selected time points; we use the $n$ chosen values to form noisy observations for asset $i$. The selection procedure is applied to $p$ assets for obtaining data on all assets. Because of random selection, the obtained data are nonsynchronized and have $n$ observations for each asset. As in the synchronized case, the true $\boldsymbol{\Gamma}$ is computed by $\sum_{r=1}^{n} \boldsymbol{\gamma}(t_{3r})/n$. But we use the generated nonsynchronized data to evaluate $\tilde{\boldsymbol{\Gamma}}$, $\mathcal{B}_b[\tilde{\boldsymbol{\Gamma}}]$ and $\mathcal{T}_\varpi[\tilde{\boldsymbol{\Gamma}}]$, where the values of $b$ and $\varpi$ are selected as before by minimizing their respective MSEs. Again we repeat the whole simulation procedure 500 times and evaluate MSEs of $\tilde{\boldsymbol{\Gamma}}$, $\mathcal{B}_b[\tilde{\boldsymbol{\Gamma}}]$ and $\mathcal{T}_\varpi[\tilde{\boldsymbol{\Gamma}}]$ based on the 500 repetitions.

5.4. *Simulation results.* In the simulations we have tried on different combinations of values for $n$ and $p$. This section displays the simulation results and reports findings based on $n = 200$ and $p = 512$. Figure 1 plots the sample paths of $\kappa(t)$ corresponding to six initial values $\kappa(0) = 0.537$, 0.762, 0.905, 0.964, 0.980, 0.995. The fixed initial values were chosen to obtain various patterns for $\boldsymbol{\Gamma}$. The figure shows that process $\kappa(t)$ is heavily influenced by its initial value, and its whole path stays within a narrow band around the initial value. Figure 2 plots the images of $\boldsymbol{\Gamma}$ corresponding to these initial values. The image plots indicate that the elements of $\boldsymbol{\Gamma}$ geometrically decay from its diagonal with rate $\kappa(t)$. The significant elements of $\boldsymbol{\Gamma}$ fall into a band along its diagonal, and its off-diagonal elements outside the band are negligible. For small $\kappa(t)$, the decay is very fast, and the band is very narrow. As $\kappa(t)$ increases, the decay gets slower and slower, and the band becomes wider and wider. As a result, $\boldsymbol{\Gamma}$ becomes less sparse and is more diffuse along its diagonal. As we will see later, it will be more difficult to estimate $\boldsymbol{\Gamma}$. Note that $\boldsymbol{\Gamma}$ is inhomogeneous, and the band is wider at lower right corner than at upper left corner. This is clearly demonstrated in Figure 2(f) for the case with $\kappa(0) = 0.995$.

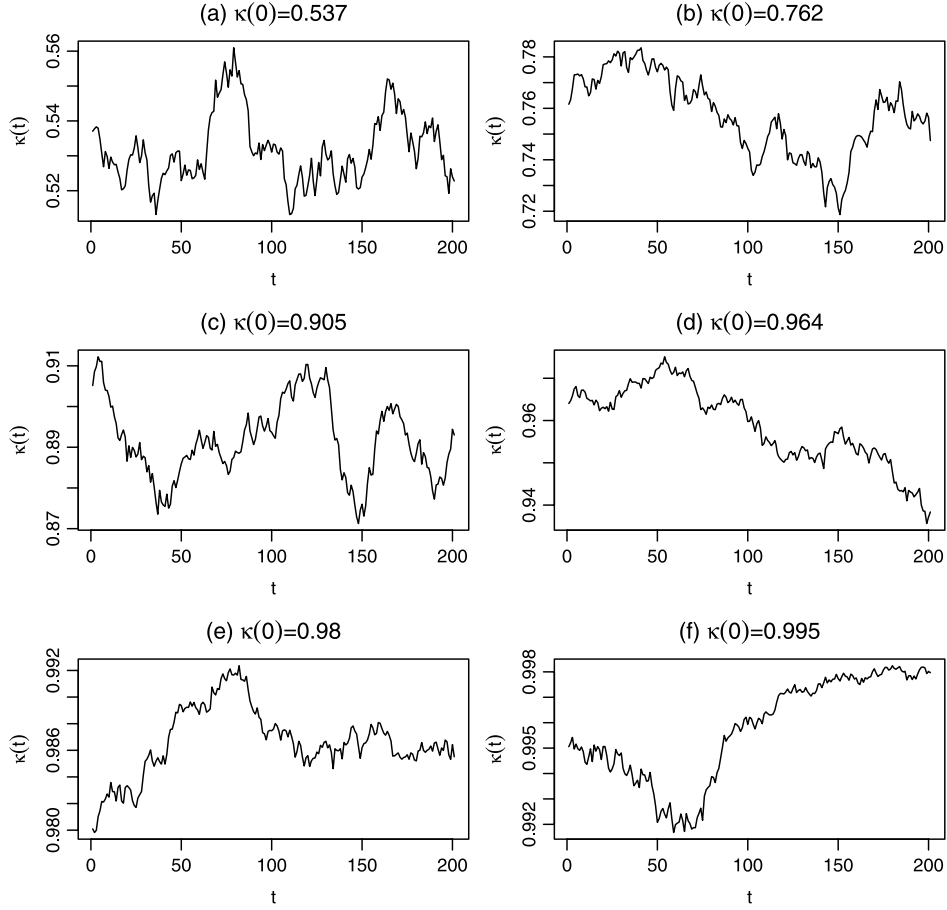

FIG. 1. *Sample path plots for the process $\kappa(t)$ under different initial values. (a)–(f) correspond to the sample paths of $\kappa(t)$ with $\kappa(0) = 0.537, 0.762, 0.905, 0.964, 0.980, 0.995$, respectively.*

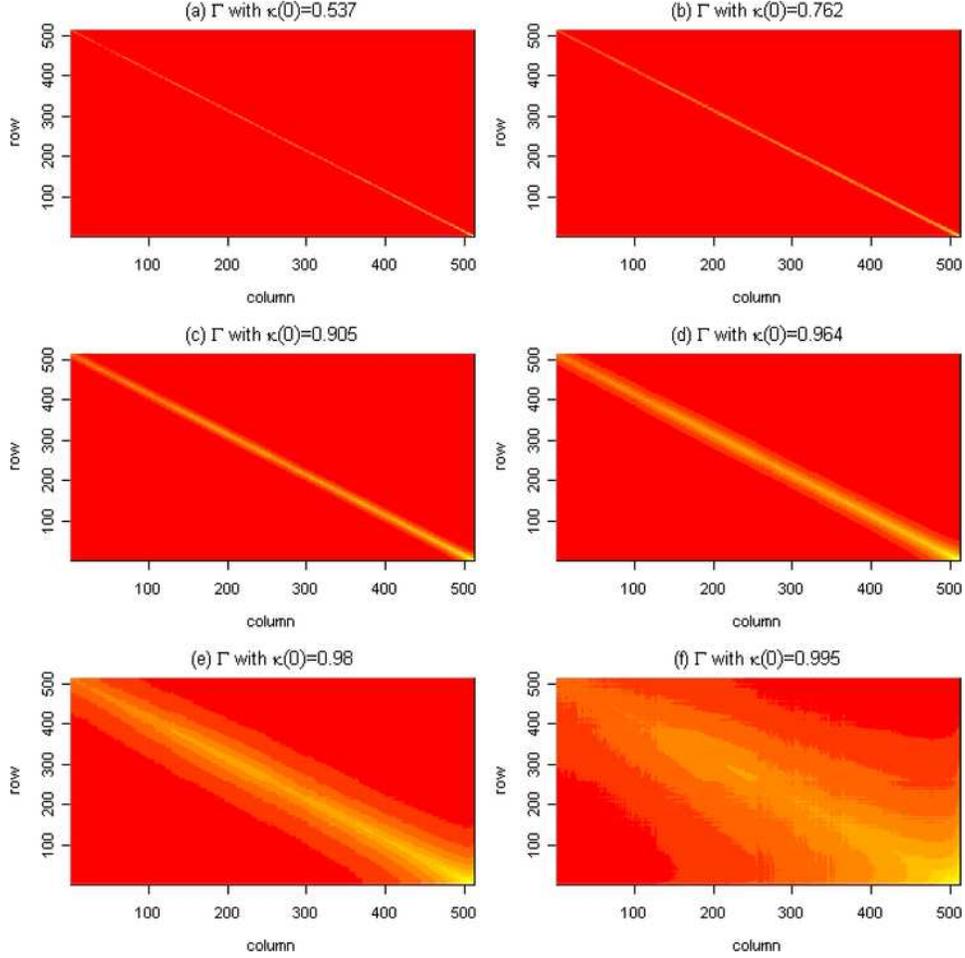

Fig. 2. *Image plots of matrix* $\Gamma$ *generated with different initial values for* $\kappa(t)$. *(a)–(f) correspond to the images of* $\Gamma$ *with* $\kappa(0) = 0.537, 0.762, 0.905, 0.964, 0.980, 0.995$, *respectively.*

The MSEs of $\tilde{\Gamma}$, $\mathcal{B}_b[\tilde{\Gamma}]$ and $\mathcal{T}_\varpi[\tilde{\Gamma}]$ are computed over combinations of six initial values of $\kappa(0)$, three noise levels and two values of $K$. The numerical results are summarized in Table 1 for the case of synchronized data.

The simulation results indicate that for the given $\Gamma$, BARVM estimator has smaller MSE than the corresponding TARVM estimator. This is due to the fact that the decay pattern of $\Gamma$ is the ideal case for banding. However, if we randomly permute the rows and columns of $\Gamma$, the resulting matrix no longer decays along its diagonal but retains the same sparsity. Such a permutation corresponds to a random shuffle of the list of stocks. For each realization matrix of $\Gamma$ displayed in Figure 2, we take a random permutation of its rows and columns. Figure 3 plots the images of the obtained matrices. The plot shows that while the significantly large elements are scattered all over the place and the decay patterns completely disappear, the sparsity remains unchanged. For such randomly permuted $\Gamma$, the TARVM estimator maintains the same performance, while the BARVM estimator performs very poorly.

TABLE 1
*MSEs of $\tilde{\Gamma}$, $\mathcal{B}_b[\tilde{\Gamma}]$ and $\mathcal{T}_\varpi[\tilde{\Gamma}]$ for noisy synchronized data*

| Noise level | Estimator | $K$ | $\kappa(0)$ | | | | | |
|---|---|---|---|---|---|---|---|---|
| | | | **0.537** | **0.762** | **0.905** | **0.964** | **0.980** | **0.995** |
| Low | $\tilde{\Gamma}$ | 1 | 5.595 | 6.039 | 7.511 | 9.959 | 12.259 | 18.270 |
| Low | $\tilde{\Gamma}$ | 5 | 10.186 | 11.80 | 14.85 | 18.79 | 21.488 | 31.85 |
| Low | $\mathcal{B}_b[\tilde{\Gamma}]$ | 1 | 0.663 | 1.195 | 3.154 | 6.000 | 8.018 | |
| Low | $\mathcal{T}_\varpi[\tilde{\Gamma}]$ | 1 | 0.845 | 2.456 | 4.595 | 7.457 | 10.928 | |
| Low | $\mathcal{B}_b[\tilde{\Gamma}]$ | 5 | 1.085 | 2.008 | 5.075 | 10.160 | 15.299 | |
| Low | $\mathcal{T}_\varpi[\tilde{\Gamma}]$ | 5 | 1.744 | 3.077 | 6.855 | 12.683 | 18.111 | |
| Medium | $\tilde{\Gamma}$ | 1 | 5.641 | 6.097 | 7.649 | 10.479 | 12.280 | 18.398 |
| Medium | $\tilde{\Gamma}$ | 5 | 10.229 | 11.81 | 14.74 | 19.17 | 21.851 | 31.95 |
| Medium | $\mathcal{B}_b[\tilde{\Gamma}]$ | 1 | 0.694 | 1.224 | 2.785 | 6.442 | 9.059 | |
| Medium | $\mathcal{T}_\varpi[\tilde{\Gamma}]$ | 1 | 0.871 | 2.466 | 4.101 | 7.680 | 12.058 | |
| Medium | $\mathcal{B}_b[\tilde{\Gamma}]$ | 5 | 1.093 | 2.022 | 4.499 | 10.302 | 15.384 | |
| Medium | $\mathcal{T}_\varpi[\tilde{\Gamma}]$ | 5 | 1.757 | 3.083 | 7.014 | 12.742 | 19.083 | |
| High | $\tilde{\Gamma}$ | 1 | 5.769 | 6.234 | 7.717 | 10.521 | 12.942 | 19.26 |
| High | $\tilde{\Gamma}$ | 5 | 10.271 | 11.86 | 14.89 | 18.85 | 21.968 | 33.54 |
| High | $\mathcal{B}_b[\tilde{\Gamma}]$ | 1 | 0.723 | 1.258 | 3.105 | 6.298 | 10.094 | |
| High | $\mathcal{T}_\varpi[\tilde{\Gamma}]$ | 1 | 0.896 | 2.429 | 4.043 | 8.765 | 12.844 | |
| High | $\mathcal{B}_b[\tilde{\Gamma}]$ | 5 | 1.077 | 2.125 | 4.601 | 10.042 | 16.041 | |
| High | $\mathcal{T}_\varpi[\tilde{\Gamma}]$ | 5 | 1.628 | 3.101 | 6.943 | 12.998 | 19.194 | |

The simulation results show that the MSEs increase in $\kappa(0)$. This can be explained by the fact that as $\kappa(0)$ increases, $\Gamma$ decays more slowly and becomes less sparse, and thus it is more difficult to estimate $\Gamma$. In fact, for $\kappa(0) = 0.995$, $\Gamma$ is so diffuse that banding and thresholding result in almost no reduction in MSEs. In other words, because $\Gamma$ is not even nearly sparse, when applying banding and thresholding procedures to $\tilde{\Gamma}$, we select almost all elements in $\tilde{\Gamma}$, and the resulting BARVM and TARVM estimators are basically the same as $\tilde{\Gamma}$. Also it is interesting to see that the MSEs increases much faster in $\kappa(0)$ than in noise level. For the chosen range of $\kappa(0)$ and the specified noise levels, $\kappa(t)$ governing the sparsity has more influence on the MSEs than noise.

The simulation results suggest that for all three noise levels, the estimators with $K = 1$ have better performance than with $K = 5$. We have tried to increase noise standard deviation up to $0.6\sqrt{\theta_i}$ and found that for noise standard deviation from $0.4\sqrt{\theta_i}$ on, the estimators with $K = 1$ perform worse than with $K = 5$. We have also tried on different values for $K$ and obtained the similar results. Like the TSRV estimator in Zhang, Mykland and Aït-Sahalia (2005), the role of $K$ is to balance subsampling and averaging in $\tilde{\Gamma}$ for the purpose of noise reduction. Its effect is clearly demonstrated by asymptotic analysis as illustrated in Theorem 1 and Remark 1. However, the simulation results imply that it needs large noise to manifest numerically the benefit of choosing $K$ greater than 1 in the construction of $\tilde{\Gamma}$.

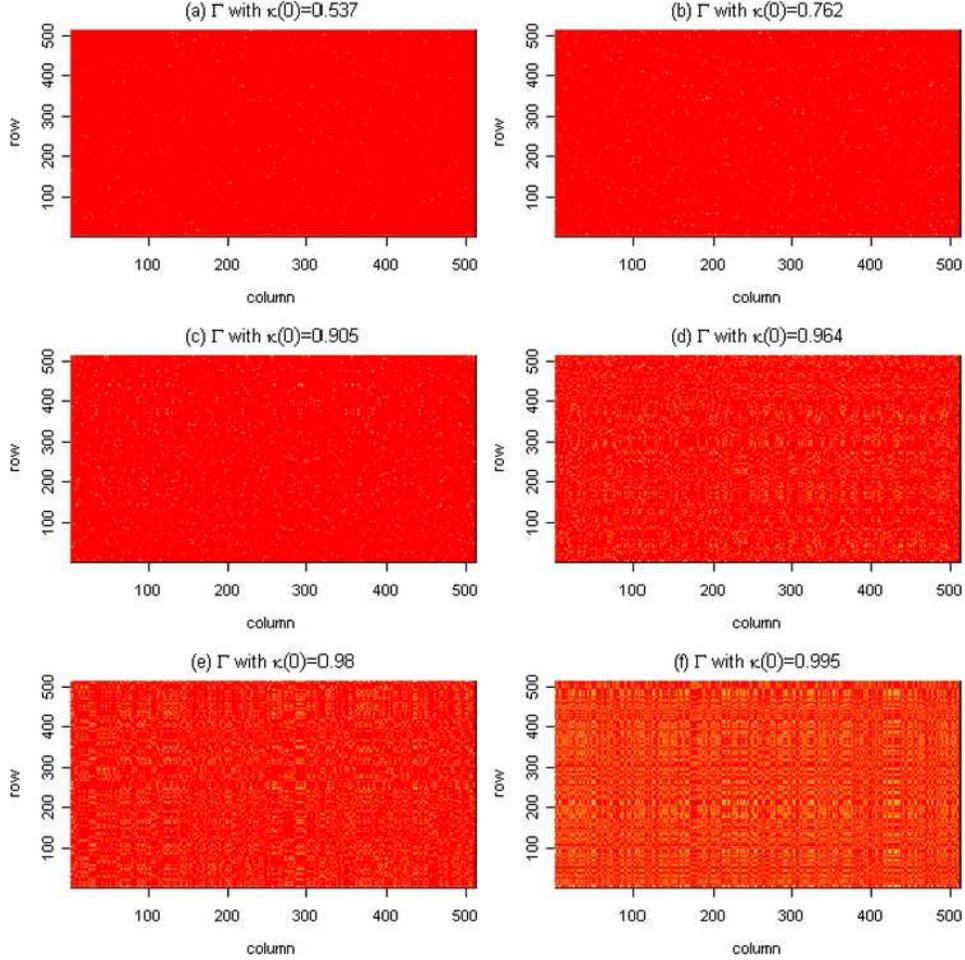

Fig. 3. *Image plots of the matrices obtained by randomly permuting rows and columns of $\Gamma$ in Figure 2. (a)–(f) correspond to the images of randomly permuted $\Gamma$ with $\kappa(0) = 0.537, 0.762, 0.905, 0.964, 0.980, 0.995$, respectively.*

TABLE 2
*MSEs of $\tilde{\boldsymbol{\Gamma}}$, $\mathcal{B}_b[\tilde{\boldsymbol{\Gamma}}]$ and $\mathcal{T}_\varpi[\tilde{\boldsymbol{\Gamma}}]$ for noisy nonsynchronized data*

| Noise level | Estimator | $K$ | $\kappa(0)$ | | | | | |
|---|---|---|---|---|---|---|---|---|
| | | | **0.537** | **0.762** | **0.905** | **0.964** | **0.980** | **0.995** |
| Low | $\tilde{\boldsymbol{\Gamma}}$ | 1 | 12.86 | 16.99 | 27.13 | 48.37 | 68.77 | 152.75 |
| Low | $\mathcal{B}_b[\tilde{\boldsymbol{\Gamma}}]$ | 1 | 3.305 | 4.778 | 8.718 | 20.375 | 34.29 | 85.10 |
| Low | $\mathcal{T}_\varpi[\tilde{\boldsymbol{\Gamma}}]$ | 1 | 3.842 | 5.281 | 13.38 | 30.29 | 46.11 | 121.15 |
| Medium | $\tilde{\boldsymbol{\Gamma}}$ | 1 | 12.98 | 17.10 | 27.15 | 48.57 | 69.23 | 153.57 |
| Medium | $\mathcal{B}_b[\tilde{\boldsymbol{\Gamma}}]$ | 1 | 3.911 | 4.834 | 8.759 | 22.49 | 36.04 | 92.70 |
| Medium | $\mathcal{T}_\varpi[\tilde{\boldsymbol{\Gamma}}]$ | 1 | 3.374 | 4.728 | 11.662 | 30.53 | 50.64 | 123.22 |
| High | $\tilde{\boldsymbol{\Gamma}}$ | 1 | 13.16 | 17.15 | 27.50 | 48.07 | 71.44 | 151.85 |
| High | $\mathcal{B}_b[\tilde{\boldsymbol{\Gamma}}]$ | 1 | 3.443 | 4.776 | 8.779 | 21.946 | 42.23 | 69.27 |
| High | $\mathcal{T}_\varpi[\tilde{\boldsymbol{\Gamma}}]$ | 1 | 3.997 | 4.902 | 11.70 | 29.98 | 56.27 | 100.13 |

Table 2 displays the MSEs of $\tilde{\boldsymbol{\Gamma}}$, $\mathcal{B}_b[\tilde{\boldsymbol{\Gamma}}]$ and $\mathcal{T}_\varpi[\tilde{\boldsymbol{\Gamma}}]$ for noisy nonsynchronized data. The comparison of Tables 1 and 2 shows that the MSEs in Table 2 are much larger than the corresponding ones in Table 1 for all three noise levels and six values of $\kappa(0)$ considered. The phenomenon suggests that the contribution in MSEs due to nonsynchronization dominates over that due to noise. We note particularly that even for the case of $\kappa(0) = 0.995$ where $\boldsymbol{\Gamma}$ is very diffuse and regularizations have little improvement on $\tilde{\boldsymbol{\Gamma}}$ for the synchronization case, regularizations in the non-synchronization case still improve $\tilde{\boldsymbol{\Gamma}}$ with sizable reduction on MSEs. Similar to the synchronized case, the estimators with $K = 1$ have smaller MSEs than with $K = 5$ for all three noise levels. Since nonsynchronization significantly inflates MSEs, it requires very large noise to manifest numerically the effect on MSE reduction by using $K > 1$ in the construction of $\tilde{\boldsymbol{\Gamma}}$. Apart from the phenomenon due to nonsynchronization, Table 2 exhibits the similar MSE patterns as Table 1.

5.5. *An application to the stock data.* We applied the proposed method to the high-frequency stock price data from the Shanghai market. Denote by $\tilde{\Gamma}_i$, $i = 1, \ldots, 177$, the daily ARVM estimators obtained in Section 5.1. For each of $\tilde{\Gamma}_i$, we computed its eigenvalues and collected them as a set. The eigenvalue sets for $\tilde{\Gamma}_i$ have wide ranges, with some very big positive values for the largest eigenvalues and many negative values for the smallest eigenvalues. As stocks have no natural ordering, the decay assumption is not realistic for volatility matrices, and banding may not be appropriate for $\tilde{\Gamma}_i$. We regularized $\tilde{\Gamma}_i$ by thresholding. The threshold $\varpi_i$ applied to $\tilde{\Gamma}_i$ was calibrated through the quantiles of the absolute entries of $\tilde{\Gamma}_i$. For $a \in (0, 1)$, let $\varpi_{i,a}$ be the $a$-quantile of the absolute entries of $\tilde{\Gamma}_i$. Then we reduced threshold selection to the selection of $a$. Define

$$(24) \qquad \Lambda(a) = \sum_{i=1}^{176} \|\tilde{\Gamma}_{i+1} - \mathcal{T}_{\varpi_{i,a}}[\tilde{\Gamma}_i]\|_2^2.$$

We selected the value of $a$ by minimizing $\Lambda(a)$ over $a \in (0, 1)$. The threshold selection procedure was motivated as follows. Because the ARVM estimators were evaluated at the daily level where stationarity is a reasonable assumption on volatility in financial time series, we predicted one day ahead of the daily realized volatility matrix by the current thresholded daily realized volatility matrix with prediction performance measured by the $\ell_2$-norm of the prediction error. Thresholds $\varpi_{i,a}$ were then selected by minimizing $\Lambda(a)$, the sum of the squared $\ell_2$-norms of the prediction errors over 176 pairs of consecutive days. Our calculation resulted in selecting 0.95 for $a$. We thresholded $\tilde{\Gamma}_i$ by $\varpi_{i,0.95}$, namely, for each of $\tilde{\Gamma}_i$ we retained its top 5% entries in magnitude and replaced the rest by zero. We evaluated the eigenvalues of $\mathcal{T}_{\varpi_{i,0.95}}[\tilde{\Gamma}_i]$. Thresholding attributes to significantly narrowing down the eigenvalue ranges. Since many $\tilde{\Gamma}_i$ had negative smallest eigenvalues, we truncated the negative eigenvalues at zero and plotted in Figure [4] the corresponding largest eigenvalues of $\tilde{\Gamma}_i$ and $\mathcal{T}_{\varpi_{i,0.95}}[\tilde{\Gamma}_i]$. The plot shows that the reductions of the largest eigenvalues due to thresholding are over 50% for many days. Since eigen based analyses like clustering analysis, principal component analysis and portfolio allocation are routinely applied in practice, the study indicates that blind applications of such analyses to large realized volatility matrices without regularization may end up with very misleading conclusions.

**6. Proofs of Theorems [2] and [3].** Denote by $C$ a generic constant whose value is free of $n$ and $p$ and may change from appearance to appearance. $O_P$ and $o_P$ denote orders in probability as both $n$ and $p$ go to infinity.

We will prove Theorem [1] in Section [7]. This section assumes (11) in Theorem [1]. We use it to establish Theorems [2] and [3] by following Bickel and Levina ([2008a], [2008b]).

PROOF OF THEOREM [3]. Using the relationship between $\ell_2$- and $\ell_\infty$-norms, triangle inequality, and decaying condition (10), we have

$$(25) \qquad \begin{aligned} \|\mathcal{B}_b[\tilde{\Gamma}] - \Gamma\|_2 &\le \|\mathcal{B}_b[\tilde{\Gamma}] - \Gamma\|_\infty \\ &\le \|\mathcal{B}_b[\tilde{\Gamma}] - \mathcal{B}_b[\Gamma]\|_\infty + \|\mathcal{B}_b[\Gamma] - \Gamma\|_\infty, \end{aligned}$$

$$(26) \qquad \|\mathcal{B}_b[\Gamma] - \Gamma\|_\infty \le \max_{1 \le i \le p} \sum_{|i-j| > b} |\Gamma_{ij}| \le 2M \sum_{k=b+1}^{\infty} k^{-\alpha-1} \le \frac{2M}{\alpha} b^{-\alpha},$$

$$(27) \qquad \|\mathcal{B}_b[\tilde{\Gamma}] - \mathcal{B}_b[\Gamma]\|_\infty \le (2b+1) \max_{|i-j| \le b} |\tilde{\Gamma}_{ij} - \Gamma_{ij}|.$$

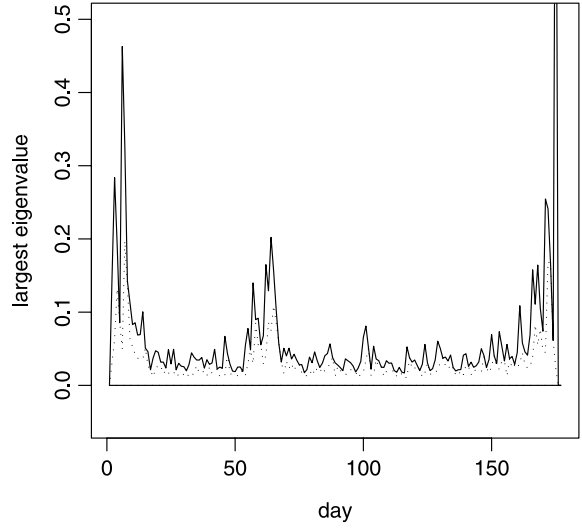

FIG. 4. *Plots of the largest eigenvalues of daily realized volatility matrices and the thresholded daily realized volatility matrices for the stock data from the Shanghai market. The solid line represents the largest eigenvalues of daily realized volatility matrices, and the dotted line corresponds to the largest eigenvalues of the thresholded daily realized volatility matrices.*

From (11) we get

$$P\Big(\max_{|i-j|\le b}|\tilde\Gamma_{ij}-\Gamma_{ij}|>d\Big)$$

$$\le \sum_{|i-j|\le b} P(|\tilde\Gamma_{ij}-\Gamma_{ij}|>d)$$

$$\le \frac{p(2b+1)}{d^\beta}\max_{|i-j|\le b}E[|\tilde\Gamma_{ij}-\Gamma_{ij}|^\beta]\le \frac{Cpbe_n^\beta}{d^\beta}.$$

Combining above probability inequality with (27) we obtain

$$P(\|\mathcal{B}_b[\tilde{\mathbf\Gamma}]-\mathcal{B}_b[\mathbf\Gamma]\|_\infty>d)\le P\Big(\max_{|i-j|\le b}|\tilde\Gamma_{ij}-\Gamma_{ij}|>d/(2b+1)\Big)$$

$$(28)$$

$$\le \frac{Cpb^{\beta+1}e_n^\beta}{d^\beta}\le \frac{Cb^{-\alpha\beta}}{d^\beta},$$

where the last inequality is due to the fact that for the selected $b$ in the theorem,

$$pb^{\beta+1}e_n^\beta\sim b^{-\alpha\beta}.$$

Collecting together (25), (26) and (28) and taking

$$d = 2d_1 b^{-\alpha} \sim (p e_n^{\beta})^{\alpha/(\alpha\beta+\beta+1)},$$

we conclude

$$P(\|\mathcal{B}_b[\tilde{\Gamma}] - \Gamma\|_\infty > d) \leq P(\|\mathcal{B}_b[\tilde{\Gamma}] - \mathcal{B}_b[\Gamma]\|_\infty > d_1 b^{-\alpha})$$
$$+ P(\|\mathcal{B}_b[\Gamma] - \Gamma\|_\infty > d_1 b^{-\alpha})$$
$$\leq \frac{C}{d_1^\beta} + P(M > \alpha d_1/2)$$
$$\leq \frac{C}{d_1^\beta} + \frac{2E[M]}{\alpha d_1} \to 0, \qquad \text{as } d_1 \to \infty.$$

This completes the proof of Theorem 3.  □

We need the following lemmas for proving Theorem 2.

LEMMA 1.  *Assume that $\Gamma$ satisfies sparse condition (10) and $\varpi$ is chosen as in Theorem 2. Then for any fixed $a > 0$,*

$$(29) \quad \max_{1 \leq i \leq p} \sum_{j=1}^p |\Gamma_{ij}| 1(|\Gamma_{ij}| \leq a\varpi) \leq a^{1-\delta} M\pi(p)\varpi^{1-\delta} = O_P(\pi(p)\varpi^{1-\delta}),$$

$$(30) \quad \max_{1 \leq i \leq p} \sum_{j=1}^p 1(|\Gamma_{ij}| \geq a\varpi) \leq a^{-\delta} M\pi(p)\varpi^{-\delta} = O_P(\pi(p)\varpi^{-\delta}).$$

PROOF.  Simple algebraic manipulation shows

$$\max_{1 \leq i \leq p} \sum_{j=1}^p |\Gamma_{ij}| 1(|\Gamma_{ij}| \leq a\varpi)$$

$$\leq (a\varpi)^{1-\delta} \max_{1 \leq i \leq p} \sum_{j=1}^p |\Gamma_{ij}|^\delta 1(|\Gamma_{ij}| \leq a\varpi)$$

$$\leq a^{1-\delta} \varpi^{1-\delta} M\pi(p) = O_P(\pi(p)\varpi^{1-\delta})$$

which proves (29). (30) is proved as follows:

$$\max_{1 \leq i \leq p} \sum_{j=1}^p 1(|\Gamma_{ij}| \geq a\varpi) \leq \max_{1 \leq i \leq p} \sum_{j=1}^p [|\Gamma_{ij}|/(a\varpi)]^\delta 1(|\Gamma_{ij}| \geq a\varpi)$$

$$\leq (a\varpi)^{-\delta} \max_{1 \leq i \leq p} \sum_{j=1}^p |\Gamma_{ij}|^\delta$$

$$\leq (a\varpi)^{-\delta} M\pi(p) = O_P(\pi(p)\varpi^{-\delta}). \qquad \square$$

LEMMA 2. *Assume that* $\boldsymbol{\Gamma}$ *satisfies sparse condition (10),* $\varpi$ *is chosen as in Theorem 2 and (11) is held. Then*

$$\max_{1 \leq i,j \leq p} |\tilde{\Gamma}_{ij} - \Gamma_{ij}| = O_P(e_n p^{2/\beta}) = o_P(\varpi), \tag{31}$$

$$P\left(\max_{1 \leq i \leq p} \sum_{j=1}^{p} 1\{|\tilde{\Gamma}_{ij} - \Gamma_{ij}| \geq \varpi/2\} > 0\right) = o(1), \tag{32}$$

$$\max_{1 \leq i \leq p} \sum_{j=1}^{p} 1(|\tilde{\Gamma}_{ij}| \geq \varpi, |\Gamma_{ij}| < \varpi) \leq 2^{\delta} M \pi(p) \varpi^{-\delta} + o_P(1)$$

$$\tag{33}$$

$$= O_P(\pi(p)\varpi^{-\delta}).$$

PROOF. Taking $d = d_1 p^{2/\beta} e_n$, applying the Markov inequality and using (11) we obtain

$$P\left(\max_{1 \leq i,j \leq p} |\tilde{\Gamma}_{ij} - \Gamma_{ij}| > d\right) \leq \sum_{i,j=1}^{p} P(|\tilde{\Gamma}_{ij} - \Gamma_{ij}| > d) \leq \frac{Cp^2 e_n^{\beta}}{d^{\beta}} = \frac{C}{d_1^{\beta}} \to 0$$

as $n, p \to \infty$ and then $d_1 \to \infty$. This proves (31). Using above inequality we have

$$P\left(\max_{1 \leq i \leq p} \sum_{j=1}^{p} 1\{|\tilde{\Gamma}_{ij} - \Gamma_{ij}| \geq \varpi/2\} > 0\right) \leq P\left(\max_{1 \leq i,j \leq p} |\tilde{\Gamma}_{ij} - \Gamma_{ij}| \geq \varpi/2\right)$$

$$\leq \frac{2^{\beta} Cp^2 e_n^{\beta}}{\varpi^{\beta}} = \frac{2^{\beta} C}{h_{n,p}^{\beta}} \to 0$$

since $h_{n,p} \to \infty$ as $n, p \to \infty$ which proves (32). To show (33) we employ (30) and (32) to get

$$\max_{1 \leq i \leq p} \sum_{j=1}^{p} 1(|\tilde{\Gamma}_{ij}| \geq \varpi, |\Gamma_{ij}| < \varpi) \leq \max_{1 \leq i \leq p} \sum_{j=1}^{p} 1(|\tilde{\Gamma}_{ij}| \geq \varpi, |\Gamma_{ij}| \leq \varpi/2)$$

$$+ \max_{1 \leq i \leq p} \sum_{j=1}^{p} 1(|\tilde{\Gamma}_{ij}| \geq \varpi, \varpi/2 < |\Gamma_{ij}| < \varpi)$$

$$\leq \max_{1 \leq i \leq p} \sum_{j=1}^{p} 1(|\tilde{\Gamma}_{ij} - \Gamma_{ij}| \geq \varpi/2)$$

$$+ \max_{1 \leq i \leq p} \sum_{j=1}^{p} 1(|\Gamma_{ij}| > \varpi/2)$$

$$\leq o_P(1) + 2^\delta M \pi(p) \varpi^{-\delta}$$
$$= O_P(\pi(p)\varpi^{-\delta}). \qquad \square$$

PROOF OF THEOREM 2. With the relationship between $\ell_2$- and $\ell_\infty$-norms and the triangle inequality, we have

$$\|\mathcal{T}_\varpi[\tilde{\boldsymbol{\Gamma}}] - \boldsymbol{\Gamma}\|_2 \leq \|\mathcal{T}_\varpi[\tilde{\boldsymbol{\Gamma}}] - \mathcal{T}_\varpi[\boldsymbol{\Gamma}]\|_2 + \|\mathcal{T}_\varpi[\boldsymbol{\Gamma}] - \boldsymbol{\Gamma}\|_2$$
$$\leq \|\mathcal{T}_\varpi[\tilde{\boldsymbol{\Gamma}}] - \mathcal{T}_\varpi[\boldsymbol{\Gamma}]\|_\infty + \|\mathcal{T}_\varpi[\boldsymbol{\Gamma}] - \boldsymbol{\Gamma}\|_\infty.$$

Lemma 1 implies

$$\|\mathcal{T}_\varpi[\boldsymbol{\Gamma}] - \boldsymbol{\Gamma}\|_\infty = \max_{1 \leq i \leq p} \sum_{j=1}^p |\Gamma_{ij}| 1(|\Gamma_{ij}| \leq \varpi) = O_P(\pi(p)\varpi^{1-\delta}).$$

The theorem is proved if we show that $\|\mathcal{T}_\varpi[\tilde{\boldsymbol{\Gamma}}] - \mathcal{T}_\varpi[\boldsymbol{\Gamma}]\|_\infty$ is also of order $\varpi^{1-\delta}\pi(p)$ in probability. Indeed, with simple algebra and Lemmas 1 and 2 we establish it as follows:

$$\|\mathcal{T}_\varpi[\tilde{\boldsymbol{\Gamma}}] - \mathcal{T}_\varpi[\boldsymbol{\Gamma}]\|_\infty \leq \max_{1 \leq i \leq p} \sum_{j=1}^p |\tilde{\Gamma}_{ij} - \Gamma_{ij}| 1(|\tilde{\Gamma}_{ij}| \geq \varpi, |\Gamma_{ij}| \geq \varpi)$$

$$+ \max_{1 \leq i \leq p} \sum_{j=1}^p |\tilde{\Gamma}_{ij}| 1(|\tilde{\Gamma}_{ij}| \geq \varpi, |\Gamma_{ij}| < \varpi)$$

$$+ \max_{1 \leq i \leq p} \sum_{j=1}^p |\Gamma_{ij}| 1(|\tilde{\Gamma}_{ij}| < \varpi, |\Gamma_{ij}| \geq \varpi)$$

$$\leq \max_{1 \leq i,j \leq p} |\tilde{\Gamma}_{ij} - \Gamma_{ij}| \max_{1 \leq i \leq p} \sum_{j=1}^p 1(|\Gamma_{ij}| \geq \varpi)$$

$$+ \max_{1 \leq i \leq p} \sum_{j=1}^p |\Gamma_{ij}| 1(|\Gamma_{ij}| < \varpi)$$

$$+ \max_{1 \leq i,j \leq p} |\tilde{\Gamma}_{ij} - \Gamma_{ij}| \max_{1 \leq i \leq p} \sum_{j=1}^p 1(|\tilde{\Gamma}_{ij}| \geq \varpi, |\Gamma_{ij}| < \varpi)$$

$$+ \varpi \max_{1 \leq i \leq p} \sum_{j=1}^p 1(|\Gamma_{ij}| \geq \varpi)$$

$$= o_P(\varpi) O_P(\pi(p)\varpi^{-\delta}) + O_P(\pi(p)\varpi^{1-\delta})$$
$$+ o_P(\varpi) O_P(\pi(p)\varpi^{-\delta}) + \varpi O_P(\pi(p)\varpi^{-\delta})$$
$$= O_P(\pi(p)\varpi^{1-\delta}),$$

where the orders in line five of the six equation array are from (29)–(31) and (33). □

## 7. Proof of Theorem 1.

7.1. *Decomposition of $\widehat{\Gamma}$ defined in (6).* We decompose $\widehat{\Gamma} = (\widehat{\Gamma}_{ij})$ into three parts with 12 terms in this subsection and show that all parts meet condition (11) in next four subsections. Denote by

$$\mathbf{Y}_r^k = (Y_1(\tau_{1,r}^k), \ldots, Y_p(\tau_{p,r}^k))^T, \qquad \mathbf{X}_r^k = (X_1(\tau_{1,r}^k), \ldots, X_p(\tau_{p,r}^k))^T,$$

$$\boldsymbol{\varepsilon}_r^k = (\varepsilon_1(\tau_{1,r}^k), \ldots, \varepsilon_p(\tau_{p,r}^k))^T$$

the vectors formed by data, true log price and noise at time points prior to $\tau_r$ for all $p$ assets. Note that since $\tau_{1,r}^k, \ldots, \tau_{p,r}^k$ are not equal, these vectors are nonsynchronized in the sense that their coordinates are the corresponding processes evaluated at different time points. From (6), we have

$$\widehat{\Gamma} = \frac{1}{K} \sum_{k=1}^{K} \sum_{r=1}^{m} [\mathbf{Y}_r^k - \mathbf{Y}_{r-1}^k][\mathbf{Y}_r^k - \mathbf{Y}_{r-1}^k]^T$$

$$= \frac{1}{K} \sum_{k=1}^{K} \sum_{r=1}^{m} [\mathbf{X}_r^k - \mathbf{X}_{r-1}^k + \boldsymbol{\varepsilon}_r^k - \boldsymbol{\varepsilon}_{r-1}^k][\mathbf{X}_r^k - \mathbf{X}_{r-1}^k + \boldsymbol{\varepsilon}_r^k - \boldsymbol{\varepsilon}_{r-1}^k]^T$$

$$(34) \quad = \frac{1}{K} \sum_{k=1}^{K} \sum_{r=1}^{m} \{[\mathbf{X}_r^k - \mathbf{X}_{r-1}^k][\mathbf{X}_r^k - \mathbf{X}_{r-1}^k]^T + [\boldsymbol{\varepsilon}_r^k - \boldsymbol{\varepsilon}_{r-1}^k][\boldsymbol{\varepsilon}_r^k - \boldsymbol{\varepsilon}_{r-1}^k]^T$$

$$+ [\mathbf{X}_r^k - \mathbf{X}_{r-1}^k][\boldsymbol{\varepsilon}_r^k - \boldsymbol{\varepsilon}_{r-1}^k]^T + [\boldsymbol{\varepsilon}_r^k - \boldsymbol{\varepsilon}_{r-1}^k][\mathbf{X}_r^k - \mathbf{X}_{r-1}^k]^T \}$$

$$\equiv \frac{1}{K} \sum_{k=1}^{K} \sum_{r=1}^{m} [\mathbf{X}_r^k - \mathbf{X}_{r-1}^k][\mathbf{X}_r^k - \mathbf{X}_{r-1}^k]^T + \mathbf{G}(1) + \mathbf{G}(2) + \mathbf{G}(3),$$

where $\mathbf{G}(1), \mathbf{G}(2), \mathbf{G}(3)$ are sums involving with noise components and will be handled in Section 7.2, and the first term corresponds to the average realized volatility estimator based on noiseless nonsynchronized true log prices $\mathbf{X}_r^k$ which will be decomposed further below. Since $\mathbf{X}_r^k$ and $\mathbf{X}_{r-1}^k$ are evaluated at time points $\tau_{i,r}^k$ and $\tau_{i,r-1}^k$, and condition A2 indicates $\tau_{i,r-1}^k \le \tau_{r-1}^k < \tau_{i,r}^k \le \tau_r^k$, we insert synchronized true log prices $\mathbf{X}(\tau_r^k)$ and $\mathbf{X}(\tau_{r-1}^k)$ in between $\mathbf{X}_r^k$ and $\mathbf{X}_{r-1}^k$ and write

$$\mathbf{X}_r^k - \mathbf{X}_{r-1}^k = \mathbf{X}_r^k - \mathbf{X}(\tau_r^k) + \mathbf{X}(\tau_r^k) - \mathbf{X}(\tau_{r-1}^k) + \mathbf{X}(\tau_{r-1}^k) - \mathbf{X}_{r-1}^k.$$

Using the above expression to expand $(\mathbf{X}_r^k - \mathbf{X}_{r-1}^k)(\mathbf{X}_r^k - \mathbf{X}_{r-1}^k)^T$, we obtain the following decomposition of the first term on the right-hand side of (34):

$$\frac{1}{K}\sum_{k=1}^{K}\sum_{r=1}^{m}[\mathbf{X}_r^k - \mathbf{X}_{r-1}^k][\mathbf{X}_r^k - \mathbf{X}_{r-1}^k]^T$$

$$= \frac{1}{K}\sum_{k=1}^{K}\sum_{r=1}^{m}\{[\mathbf{X}(\tau_r^k) - \mathbf{X}(\tau_{r-1}^k)][\mathbf{X}(\tau_r^k) - \mathbf{X}(\tau_{r-1}^k)]^T$$

$$+ [\mathbf{X}_r^k - \mathbf{X}(\tau_r^k)][\mathbf{X}_r^k - \mathbf{X}(\tau_r^k)]^T$$

$$+ [\mathbf{X}(\tau_{r-1}^k) - \mathbf{X}_{r-1}^k][\mathbf{X}(\tau_{r-1}^k) - \mathbf{X}_{r-1}^k]^T$$

$$+ [\mathbf{X}_r^k - \mathbf{X}(\tau_r^k)][\mathbf{X}(\tau_{r-1}^k) - \mathbf{X}_{r-1}^k]^T$$

(35) $$+ [\mathbf{X}(\tau_{r-1}^k) - \mathbf{X}_{r-1}^k][\mathbf{X}_r^k - \mathbf{X}(\tau_r^k)]^T$$

$$+ [\mathbf{X}_r^k - \mathbf{X}(\tau_r^k)][\mathbf{X}(\tau_r^k) - \mathbf{X}(\tau_{r-1}^k)]^T$$

$$+ [\mathbf{X}(\tau_r^k) - \mathbf{X}(\tau_{r-1}^k)][\mathbf{X}_r^k - \mathbf{X}(\tau_r^k)]^T$$

$$+ [\mathbf{X}(\tau_{r-1}^k) - \mathbf{X}_{r-1}^k][\mathbf{X}(\tau_r^k) - \mathbf{X}(\tau_{r-1}^k)]^T$$

$$+ [\mathbf{X}(\tau_r^k) - \mathbf{X}(\tau_{r-1}^k)][\mathbf{X}(\tau_{r-1}^k) - \mathbf{X}_{r-1}^k]^T\}$$

$$\equiv \mathbf{V} + \mathbf{H}(1) + \cdots + \mathbf{H}(8),$$

where $\mathbf{V}$ corresponds to the average realized volatility estimator based on synchronized true log prices $\mathbf{X}(\tau_r^k)$, and $\mathbf{H}(1),\dots,\mathbf{H}(8)$ are contributions due to nonsynchronization in true log prices. Then from (34) and (35) we decompose $\tilde{\mathbf{\Gamma}} - \mathbf{\Gamma} = \hat{\mathbf{\Gamma}} - 2m\hat{\mathbf{\eta}} - \mathbf{\Gamma}$ into three parts with 12 terms,

(36) $$\tilde{\mathbf{\Gamma}} - \mathbf{\Gamma} = [\mathbf{G}(1) - 2m\hat{\mathbf{\eta}} + \mathbf{G}(2) + \mathbf{G}(3)] + [\mathbf{V} - \mathbf{\Gamma}]$$

$$+ [\mathbf{H}(1) + \cdots + \mathbf{H}(8)].$$

Propositions 1–3 in Sections 7.2–7.4 below, respectively, establish orders for the $\beta$th moments of the three parts on the right-hand side of (36). Putting these order results together and applying the Hölder inequality, we immediately prove Theorem 1.

7.2. *Analyze* $\mathbf{G}s$ *for the effect of micro-structure noise.* Let

$$\mathbf{G}(1) = (G_{ij}(1)), \qquad \mathbf{G}(2) = (G_{ij}(2)), \qquad \mathbf{G}(3) = (G_{ij}(3)).$$

The purpose of this subsection is to show

PROPOSITION 1. *Under the assumptions of Theorem 1, we have for* $1 \leq i, j \leq p$,

$$E[|G_{ij}(1) - 2m\hat{\eta}_i 1(i=j) + G_{ij}(2) + G_{ij}(3)|^\beta] \leq C[(Kn^{-1/2})^{-\beta} + K^{-\beta/2}].$$

We prove the proposition by deriving the orders for the $\beta$th moments of $G_{ij}(1) - 2m\eta_i 1(i = j)$, $G_{ij}(2)$, $G_{ij}(3)$ and $2m(\widehat{\eta}_i - \eta_i)$ in Lemmas 3–5 below.

LEMMA 3. *Under the assumptions of Theorem 1, we have for $1 \le i, j \le p$,*

$$E[|G_{ij}(1) - 2m\eta_i 1(i = j)|^\beta] \le C(Kn^{-1/2})^{-\beta}.$$

PROOF. From the definition of $\mathbf{G}(1)$ in (34), we have

$$G_{ij}(1) - 2m\eta_i 1(i = j)$$

$$= \frac{1}{K} \sum_{k=1}^{K} \sum_{r=1}^{m} [\varepsilon_i(\tau_{i,r}^k) - \varepsilon_i(\tau_{i,r-1}^k)][\varepsilon_j(\tau_{j,r}^k) - \varepsilon_j(\tau_{j,r-1}^k)] - 2\eta_i 1(i = j)$$

$$= \frac{1}{K} \sum_{k=1}^{K} \sum_{r=1}^{m} \{\varepsilon_i(\tau_{i,r}^k)\varepsilon_j(\tau_{j,r}^k) - \eta_i 1(i = j)$$

$$\qquad\qquad + \varepsilon_i(\tau_{i,r-1}^k)\varepsilon_j(\tau_{j,r-1}^k) - \eta_i 1(i = j)$$

$$\qquad\qquad - \varepsilon_i(\tau_{i,r}^k)\varepsilon_j(\tau_{j,r-1}^k) - \varepsilon_i(\tau_{i,r-1}^k)\varepsilon_j(\tau_{j,r}^k)\}$$

$$\equiv \frac{1}{K}[R_1 + R_2 - R_3 - R_4].$$

Note that $\tau_{i,r}^k$ and $\tau_{j,r}^k$ are equal to some $t_{i\ell}$ and $t_{j\ell}$, for fixed $i$, $\varepsilon_i(t_{i\ell})$ are i.i.d., and for $i \ne j$, $\{\varepsilon_i(t_{i\ell})\}$ and $\{\varepsilon_j(t_{j\ell})\}$ are independent. Thus, $R_1$ and $R_2$ are the sums of $\varepsilon_i(t_{i\ell})\varepsilon_j(t_{j\ell})$, $R_3$ is the sum of $\varepsilon_i(t_{i\ell})$ multiplying by the lagged $\varepsilon_j(t_{j\ell})$, and $R_4$ is the sum of $\varepsilon_j(t_{j\ell})$ multiplying by the lagged $\varepsilon_i(t_{i\ell})$. As a result, all four $R$s are martingales. We apply the Burkholder inequality [Chow and Teicher (1997), Section 11.2] to $R$s and use the moment condition on $\varepsilon_i(t_{i\ell})$ in condition A1 to obtain

$$E[|G_{ij}(1) - 2m\eta_i 1(i = j)|^\beta]$$

$$\le CK^{-\beta}(Km)^{\beta/2-1} \sum_{k=1}^{K} \sum_{r=1}^{m} E\{|\varepsilon_i(\tau_{i,r}^k)\varepsilon_j(\tau_{j,r}^k) - \eta_i 1(i = j)|^\beta$$

$$\qquad\qquad + |\varepsilon_i(\tau_{i,r-1}^k)\varepsilon_j(\tau_{j,r-1}^k) - \eta_i 1(i = j)|^\beta$$

$$\qquad\qquad + |\varepsilon_i(\tau_{i,r})\varepsilon_j(\tau_{j,r-1})|^\beta$$

$$\qquad\qquad\qquad + |\varepsilon_i(\tau_{i,r-1})\varepsilon_j(\tau_{j,r})|^\beta\}$$

$$\le CK^{-\beta}(Km)^{\beta/2}\{E[|\varepsilon_i(t_{i,1})\varepsilon_j(t_{j,1}) - \eta_i 1(i = j)|^\beta]$$

$$\qquad\qquad + E[|\varepsilon_i(t_{i,1})|^\beta]E[|\varepsilon_j(t_{j,1})|^\beta]\}$$

$$\le C(m/K)^{\beta/2} \le C(Kn^{-1/2})^{-\beta}$$

which proves the lemma.  □

LEMMA 4.  *Under the assumptions of Theorem* 1, *we have for* $1 \leq i, j \leq p$,

$$E[|G_{ij}(2)|^\beta] \leq CK^{-\beta/2}, \qquad E[|G_{ij}(3)|^\beta] \leq CK^{-\beta/2}.$$

PROOF.  Because of similarity, we provide arguments only for the first result. Simple algebra shows that

$$
\begin{aligned}
KG_{ij}(2) &= \sum_{k=1}^{K} \sum_{r=1}^{m} [X_i(\tau_{i,r}^k) - X_i(\tau_{i,r-1}^k)][\varepsilon_j(\tau_{j,r}^k) - \varepsilon_j(\tau_{j,r-1}^k)] \\
&= \sum_{k=1}^{K} \sum_{r=1}^{m} [X_i(\tau_{i,r}^k) - X_i(\tau_{i,r-1}^k)]\varepsilon_j(\tau_{j,r}^k) \\
&\quad - \sum_{k=1}^{K} \sum_{r=1}^{m} [X_i(\tau_{i,r}^k) - X_i(\tau_{i,r-1}^k)]\varepsilon_j(\tau_{j,r-1}^k) \\
&= \sum_{k=1}^{K} \sum_{r=1}^{m-1} [2X_i(\tau_{i,r}^k) - X_i(\tau_{i,r-1}^k) - X_i(\tau_{i,r+1}^k)]\varepsilon_j(\tau_{j,r}^k) \\
&\quad + \sum_{k=1}^{K} [X_i(\tau_{i,m}^k) - X_i(\tau_{i,m-1}^k)]\varepsilon_j(\tau_{j,m}^k) \\
&\quad - \sum_{k=1}^{K} [X_i(\tau_{i,1}^k) - X_i(\tau_{i,0}^k)]\varepsilon_j(\tau_{j,0}^k) \\
&\equiv R_5 + R_6 - R_7.
\end{aligned}
$$

(37)

Conditional on the whole path of $X_t$, $R_5$, $R_6$ and $R_7$, all are sums of independent random variables $\varepsilon_j(t_{j\ell})$. Thus

$$
\begin{aligned}
E[|R_5|^\beta | X] &\leq C(Km)^{\beta/2-1} \\
&\quad \times \sum_{k=1}^{K} \sum_{r=1}^{m-1} |2X_i(\tau_{i,r}^k) - X_i(\tau_{i,r-1}^k) - X_i(\tau_{i,r+1}^k)|^\beta E[|\varepsilon_j(\tau_{j,r}^k)|^\beta] \\
&\leq C(Km)^{\beta/2-1} \sum_{k=1}^{K} \sum_{r=1}^{m-1} |2X_i(\tau_{i,r}^k) - X_i(\tau_{i,r-1}^k) - X_i(\tau_{i,r+1}^k)|^\beta.
\end{aligned}
$$

Taking expectation in above inequality we get

$$E[|R_5|^\beta] \leq C(Km)^{\beta/2-1} \sum_{k=1}^{K} \sum_{r=1}^{m-1} E|2X_i(\tau_{i,r}^k) - X_i(\tau_{i,r-1}^k) - X_i(\tau_{i,r+1}^k)|^\beta$$

$$\leq C(Km)^{\beta/2-1} \sum_{k=1}^{K} \sum_{r=1}^{m-1} E \left| \int_{\tau_{i,r-1}^k}^{\tau_{i,r}^k} \boldsymbol{\sigma}_i(s) \, d\mathbf{B}_s - \int_{\tau_{i,r}^k}^{\tau_{i,r+1}^k} \boldsymbol{\sigma}_i(s) \, d\mathbf{B}_s \right|^{\beta}$$

$$(38) \qquad \leq C(Km)^{\beta/2-1} \sum_{k=1}^{K} \sum_{r=1}^{m-1} E \left| \int_{\tau_{i,r-1}^k}^{\tau_{i,r}^k} \gamma_{ii}(s) \, ds + \int_{\tau_{i,r}^k}^{\tau_{i,r+1}^k} \gamma_{ii}(s) \, ds \right|^{\beta/2}$$

$$\leq C(Km)^{\beta/2-1} \sum_{k=1}^{K} \sum_{r=1}^{m-1} m^{-\beta/2} \max_{0 \leq s \leq 1} E |\gamma_{ii}(s)|^{\beta/2}$$

$$\leq C(Km)^{\beta/2} m^{-\beta/2} = C K^{\beta/2},$$

where the third inequality is due to an application of the Burkholder inequality [He, Wang and Yan (1992), Section 10.5 and Jacod and Shiryaev (2003), Section 7.3] to the stochastic integrals.

Similarly, we have

$$E[|R_6|^{\beta}|X] \leq C K^{\beta/2-1} \sum_{k=1}^{K} |X_i(\tau_{i,m}^k) - X_i(\tau_{i,m-1}^k)|^{\beta} E[|\varepsilon_j(\tau_{j,m}^k)|^{\beta}]$$

$$\leq C K^{\beta/2-1} \sum_{k=1}^{K} |X_i(\tau_{i,m}^k) - X_i(\tau_{i,m-1}^k)|^{\beta},$$

$$(39) \qquad E[|R_6|^{\beta}] \leq C K^{\beta/2-1} \sum_{k=1}^{K} E |X_i(\tau_{i,m}^k) - X_i(\tau_{i,m-1}^k)|^{\beta} \leq C K^{\beta/2} m^{-\beta/2},$$

$$E[|R_7|^{\beta}|X] \leq C K^{\beta/2-1} \sum_{k=1}^{K} |X_i(\tau_{i,1}^k) - X_i(\tau_{i,0}^k)|^{\beta} E[|\varepsilon_j(\tau_{j,0}^k)|^{\beta}]$$

$$\leq C K^{\beta/2-1} \sum_{k=1}^{K} |X_i(\tau_{i,1}^k) - X_i(\tau_{i,0}^k)|^{\beta},$$

$$(40) \qquad E[|R_7|^{\beta}] \leq C K^{\beta/2} m^{-\beta/2}.$$

Collecting (37)–(40) together, we prove the result for $G_{ij}(2)$. $\square$

LEMMA 5. *Under the assumptions of Theorem 1, we have for $1 \leq i \leq p$,*

$$E[|m(\widehat{\eta}_i - \eta_i)|^{\beta}] \leq C(K n^{-1/2})^{-\beta}.$$

PROOF. Taking $K = 1$ in the proofs of Lemmas 3 and 4, we have that conditional on the whole path of $X_i(t)$,

$$E[|\widehat{\eta}_i - \eta_i|^{\beta}|X_i] \leq C n_i^{-\beta} \left( \sum_{\ell=1}^{n_i} [X_i(t_{i\ell}) - X_i(t_{i,\ell-1})]^{2\beta} + n_i^{\beta/2} \right),$$

$$E[|\widehat{\eta}_i - \eta_i|^\beta] \le Cn_i^{-\beta}\left(\sum_{\ell=1}^{n_i} E\{[X_i(t_{i\ell}) - X_i(t_{i,\ell-1})]^{2\beta}\} + n_i^{\beta/2}\right)$$

$$\le Cn_i^{-\beta}\left(\sum_{\ell=1}^{n_i} Cn_i^{-\beta} + n_i^{\beta/2}\right) \le Cn^{-\beta/2}$$

which immediately shows the lemma as $n = mK$. $\quad\square$

### 7.3. *Analyze* $\mathbf{V}$ *for average realized volatility based on synchronized true log price.* Let

$$[X_i, X_j]^{(k)} = \sum_{r=1}^m [X_i(\tau_r^k) - X_i(\tau_{r-1}^k)][X_j(\tau_r^k) - X_j(\tau_{r-1}^k)],$$

(41)

$$[\mathbf{X}, \mathbf{X}]^{(k)} = ([X_i, X_j]^{(k)}),$$

where $[X_i, X_j]^{(k)}$ is realized co-volatility between $X_i(t)$ and $X_j(t)$ based on true log prices at the same grid times $\tau_r^k$, $r = 1, \ldots, m$, and $\mathbf{V}$ is equal to the average of $K$ realized volatility matrices based on synchronized true log prices $\mathbf{X}(\tau_r^k)$, $r = 1, \ldots, m$. Then from (35), we have that

(42) $$\mathbf{V} = (V_{ij}) = \frac{1}{K}\sum_{k=1}^K [\mathbf{X}, \mathbf{X}]^{(k)}.$$

PROPOSITION 2. *Under the assumptions of Theorem* 1, *we have for* $1 \le i, j \le p$,

$$E(|V_{ij} - \Gamma_{ij}|^\beta) \le Cm^{-\beta/2}.$$

PROOF. Note that

$$V_{ij} = \frac{1}{K}\sum_{k=1}^K [X_i, X_j]^{(k)},$$

$$|V_{ij} - \Gamma_{ij}|^\beta = \left|\frac{1}{K}\sum_{k=1}^K ([X_i, X_j]^{(k)} - \Gamma_{ij})\right|^\beta$$

$$\le \frac{1}{K}\sum_{k=1}^K |[X_i, X_j]^{(k)} - \Gamma_{ij}|^\beta.$$

Lemma 6 below gives the moment inequality for each $[X_i, X_j]^{(k)}$, from which we immediately show

$$E(|V_{ij} - \Gamma_{ij}|^\beta) \le \frac{1}{K}\sum_{k=1}^K E|[X_i, X_j]^{(k)} - \Gamma_{ij}|^\beta$$

$$\leq Cm^{-\beta/2}. \qquad \square$$

LEMMA 6. *Under the assumptions of Theorem 1, we have that for $1 \leq k \leq K$ and $1 \leq i, j \leq p$,*

$$E(|[X_i, X_j]^{(k)} - \Gamma_{ij}|^\beta) \leq Cm^{-\beta/2}.$$

PROOF. Note that $[\mathbf{X}, \mathbf{X}]^{(k)}$ is defined by $\mathbf{X}(\tau_r^k)$ for $r = 1, \ldots, m$ while $k$ is fixed. For fixed $k$, $\tau_r^k$, $r = 1, \ldots, m$, are equally spaced grids on $[0, 1]$. As $k$ is fixed throughout the proof of the lemma, to simplify notation we write $\tau_r^k$ as $\tau_r$, $r = 1, \ldots, m$, in the rest of the proof. The effects of diffusion drifts on realized volatility are in high and negligible orders. For simplicity, we assume $\boldsymbol{\mu_t} = 0$. Observe that

$$
\begin{aligned}
&[X_i, X_j]^{(k)} - \Gamma_{ij} \\
(43) \quad &= \sum_{r=1}^m \{X_i(\tau_r) - X_i(\tau_{r-1})\}\{X_j(\tau_r) - X_j(\tau_{r-1})\} - \int_0^1 \gamma_{ij}(t)\, dt \\
&= \sum_{r=1}^m \left\{ \int_{\tau_{r-1}}^{\tau_r} (\boldsymbol{\sigma}_{is})^T\, d\mathbf{B}_s \int_{\tau_{r-1}}^{\tau_r} (\boldsymbol{\sigma}_{js})^T\, d\mathbf{B}_s - \int_{\tau_{r-1}}^{\tau_r} \gamma_{ij}(t)\, dt \right\} \\
&= \sum_{r=1}^m \{D_{ir}(\tau_r) D_{jr}(\tau_r) - [D_{ir}, D_{jr}]_{\tau_r}\},
\end{aligned}
$$

where $\boldsymbol{\sigma}_{is} = (\sigma_{1i,s}, \ldots, \sigma_{pi,s})^T$ is the $i$th column of $\boldsymbol{\sigma}_s$,

$$\gamma_{ij}(s) = (\boldsymbol{\sigma}_{is})^T \boldsymbol{\sigma}_{js} = \sum_{a=1}^p \sigma_{ai,s} \sigma_{aj,s},$$

$$D_{ir}(t) = X_i(\tau_r) - X_i(\tau_{r-1}) = \int_{\tau_{r-1}}^t (\boldsymbol{\sigma}_{is})^T\, d\mathbf{B}_s, \qquad t \in [\tau_{r-1}, \tau_r)$$

and $D_i(t) = D_{ir}(t)$ for $t \in [\tau_{r-1}, \tau_r)$. Applying Itô's lemma, we obtain the following stochastic integral expression:

$$(44) \quad D_{ir}(\tau_r) D_{jr}(\tau_r) - [D_{ir}, D_{jr}]_{\tau_r} = \int_{\tau_{r-1}}^{\tau_r} \{D_{ir}(s)\boldsymbol{\sigma}_{is} + D_{jr}(s)\boldsymbol{\sigma}_{js}\}^T\, d\mathbf{B}_s$$

which has quadratic variation

$$(45) \quad \int_{\tau_{r-1}}^{\tau_r} \{D_{ir}^2(s)\gamma_{ii}(s) + D_{jr}^2(s)\gamma_{jj}(s) + 2D_{ir}(s)D_{jr}(s)\gamma_{ij}(s)\}\, ds.$$

Thus from (43)–(45) we immediately show that $[X_i, X_j]^{(k)} - \Gamma_{ij}$ has quadratic variation

$$
\begin{aligned}
(46) \qquad & \int_0^1 \{D_i^2(s)\gamma_{ii}(s) + D_j^2(s)\gamma_{jj}(s) + 2D_i(s)D_j(s)\gamma_{ij}(s)\}\, ds \\
& \leq 2\int_0^1 \{D_i^2(s)\gamma_{ii}(s) + D_j^2(s)\gamma_{jj}(s)\}\, ds,
\end{aligned}
$$

where the Cauchy–Schwarz inequality is employed to show that $|\gamma_{ij}(s)|^2 \leq \gamma_{ii}(s)\gamma_{jj}(s)$, and hence

$$
2|D_i(s)D_j(s)\gamma_{ij}(s)| \leq D_i^2(s)\gamma_{ii}(s) + D_j^2(s)\gamma_{jj}(s).
$$

Applying the Burkholder inequality [He, Wang and Yan (1992), Section 10.5 and Jacod and Shiryaev (2003), Section 7.3] to the stochastic integral expression of $[X_i, X_j]^{(k)} - \Gamma_{ij}$ given by (43)–(44) and using (44)–(46), we have

$$
\begin{aligned}
(47) \qquad E\{|[X_i, X_j]^{(k)} - \Gamma_{ij}|^\beta\} & \leq CE\left\{\left|\int_0^1 [D_i^2(s)\gamma_{ii}(s) + D_j^2(s)\gamma_{jj}(s)]\, ds\right|^{\beta/2}\right\} \\
& \leq C\int_0^1 E\{|D_i^2(s)\gamma_{ii}(s) + D_j^2(s)\gamma_{jj}(s)|^{\beta/2}\}\, ds.
\end{aligned}
$$

On the other hand, simple algebraic manipulations derive the inequalities

$$
\begin{aligned}
(48) \qquad |D_i^2(s)\gamma_{ii}(s) + D_j^2(s)\gamma_{jj}(s)|^{\beta/2} & \leq 2^{\beta/2-1}\{|D_i^2(s)\gamma_{ii}(s)|^{\beta/2} \\
& \qquad + |D_j^2(s)\gamma_{jj}(s)|^{\beta/2}\}, \\
(49) \qquad E\{|D_i^2(s)\gamma_{ii}(s)|^{\beta/2}\} & \leq \sqrt{E\{|D_i(s)|^{2\beta}\}E\{|\gamma_{ii}(s)|^\beta\}},
\end{aligned}
$$

and we apply the Burkholder inequality [He, Wang and Yan (1992), Section 10.5 and Jacod and Shiryaev (2003), Section 7.3] to the stochastic integral $D_i(s)$, and obtain

$$
\begin{aligned}
(50) \qquad E\{|D_i(s)|^{2\beta}\} & \leq C\max_{1 \leq r \leq m} E\left\{\left|\int_{\tau_{r-1}}^{\tau_r} \gamma_{ii}(s)\, ds\right|^\beta\right\} \\
& \leq C\max_{1 \leq r \leq m}(\tau_r - \tau_{r-1})^\beta \max_{0 \leq s \leq 1} E\{|\gamma_{ii}(s)|^\beta\} \\
& \leq Cm^{-\beta}\max_{0 \leq s \leq 1} E\{|\gamma_{ii}(s)|^\beta\}.
\end{aligned}
$$

Collecting together (47)–(50) we arrive at

$$
E\{|[X_i, X_j]^{(k)} - \Gamma_{ij}|^\beta\} \leq Cm^{-\beta/2}\max_{0 \leq s \leq 1} E\{|\gamma_{ii}(s)|^\beta + E|\gamma_{jj}(s)|^\beta\}.
$$

This completes the proof of the lemma. $\quad\square$

7.4. *Analyze* $\mathbf{H}s$ *for the effect of nonsynchronization.* From (35) we know that $\mathbf{H}(1) = (H_{ij}(1)), \ldots, \mathbf{H}(8) = (H_{ij}(8))$ are attributed to nonsynchronization in true log prices.

PROPOSITION 3. *Under the assumptions of Theorem 1, we have for* $1 \leq i, j \leq p$,

$$E\{|H_{ij}(1)|^\beta + \cdots + |H_{ij}(8)|^\beta\} \leq C\{(m/n)^\beta + n^{-\beta/2}\}.$$

PROOF. Because of similarity we provide arguments only for $H_{ij}(1)$, that is, to show

$$E\{|H_{ij}(1)|^\beta\} \leq C\{(m/n)^\beta + n^{-\beta/2}\}.$$

Since

$$H_{ij}(1) = \frac{1}{K} \sum_{k=1}^{K} H_{ij}^k(1),$$

where

$$H_{ij}^k(1) = \sum_{r=1}^{m} \int_{\tau_{i,r}^k}^{\tau_r^k} dX_i(s) \int_{\tau_{j,r}^k}^{\tau_r^k} dX_j(s),$$

we have

$$|H_{ij}(1)|^\beta \leq \frac{1}{K} \sum_{k=1}^{K} |H_{ij}^k(1)|^\beta.$$

So we need to show for $1 \leq k \leq K$,

$$E\{|H_{ij}^k(1)|^\beta\} \leq C\{(m/n)^\beta + n^{-\beta/2}\}. \tag{51}$$

Let

$$\bar{H}_{ij}^k(1) = \int_{\tau_{i,r}^k}^{\tau_r^k} dX_i(s) \int_{\tau_{j,r}^k}^{\tau_r^k} dX_j(s) - \int_{\tau_{i,r}^k \vee \tau_{j,r}^k}^{\tau_r^k} \gamma_{i,j} \, ds. \tag{52}$$

Note that

$$\begin{aligned}
&\sum_{r=1}^{m} \int_{\tau_{i,r}^k \vee \tau_{j,r}^k}^{\tau_r^k} \gamma_{ij}(s) \, ds \\
&= \sum_{r=1}^{m} \int_0^1 \gamma_{ij}(u(\tau_r^k - \tau_{i,r}^k \vee \tau_{j,r}^k) + \tau_r^k)(\tau_r^k - \tau_{i,r}^k \vee \tau_{j,r}^k) \, du.
\end{aligned} \tag{53}$$

The definition of $\tau_{i,r}^k$ and assumption A2 show

$$0 \leq h_{ij}^k = \tau_r^k - \tau_{i,r}^k \vee \tau_{j,r}^k \leq Cn^{-1}. \tag{54}$$

From the relationship between $H_{ij}^k(1)$ and $\bar{H}_{ij}^k(1)$ and using (53) and (54), we prove (51) as follows:

$$E\{|H_{ij}^k(1)|^\beta\} \le CE\left\{\left|\sum_{r=1}^m \int_0^1 n^{-1}|\gamma_{ij}(uh_{ij}^k + \tau_r^k)|\,du\right|^\beta\right\} + CE\{|\bar{H}_{ij}^k(1)|^\beta\}$$

$$\le Cn^{-\beta}m^{\beta-1}\sum_{r=1}^m \max_{0\le s\le 1} E\{|\gamma_{ij}(s)|^\beta\} + CE([\bar{H}_{ij}^k(1), \bar{H}_{ij}^k(1)]^{\beta/2})$$

$$\le C(m/n)^\beta + Cn^{-\beta/2} = C\{(m/n)^\beta + n^{-\beta/2}\},$$

where for the second inequality we have applied the Burkholder inequality [He, Wang and Yan (1992), Section 10.5 and Jacod and Shiryaev (2003), Section 7.3] to $E\{|\bar{H}_{ij}^k(1)|^\beta\}$, and the third inequality is due to Lemma 7 below. □

LEMMA 7. *Under the assumptions of Theorem 1, we have for $1 \le k \le K$ and $1 \le i, j \le p$,*

$$E([\bar{H}_{ij}^k(1), \bar{H}_{ij}^k(1)]^{\beta/2}) \le Cn^{-\beta/2}.$$

PROOF. With $h_{ij}^k$ defined in (54), the same argument for deriving (46) shows

$$[\bar{H}_{ij}^k(1), \bar{H}_{ij}^k(1)]$$

$$= \sum_{r=1}^m \int_{\tau_{i,r}^k \vee \tau_{j,r}^k}^{\tau_r^k} \{D_{i,r}^2(s)\gamma_{ii}(s) + D_{j,r}^2(s)\gamma_{jj}(s) + 2D_{i,r}(s)D_{j,r}\gamma_{ij}(s)\}\,ds$$

$$= \sum_{r=1}^m \int_0^1 h_{ij}^k\{D_{i,r}^2(uh_{ij}^k + \tau_r^k)\gamma_{ii}(uh_{ij}^k + \tau_r^k)$$

$$+ D_{j,r}^2(uh_{ij}^k + \tau_r^k)\gamma_{jj}(uh_{ij}^k + \tau_r^k)$$

$$+ 2D_{i,r}(uh_{ij}^k + \tau_r^k)D_{j,r}\gamma_{ij}(uh_{ij}^k + \tau_r^k)\}\,du$$

$$\le Cn^{-1}\int_0^1 \sum_{r=1}^m \{D_{i,r}^2(uh_{ij}^k + \tau_r^k)\gamma_{ii}(uh_{ij}^k + \tau_r^k)$$

$$+ D_{j,r}^2(uh_{ij}^k + \tau_r^k)\gamma_{jj}(uh_{ij}^k + \tau_r^k)\}\,du,$$

where the last inequality is due to (54) and the facts that $|\gamma_{ij}|^2 \le \gamma_{ii}\gamma_{jj}$, and thus

$$2|D_{i,r}(s)D_{j,r}(s)\gamma_{ij}(s)| \le D_{i,r}^2(s)\gamma_{ii}(s) + D_{j,r}^2(s)\gamma_{jj}(s).$$

Hence, taking $\beta/2$ power and then expectation on both sides of the inequality for $[\bar{H}_{ij}^k(1), \bar{H}_{ij}^k(1)]$, we prove the lemma as follows:

$$E([\bar{H}_{ij}^k(1), \bar{H}_{ij}^k(1)]^{\beta/2})$$

$$\leq Cn^{-\beta/2}m^{\beta/2-1}\sum_{r=1}^m \max_{0\leq s\leq 1} E\{D_{i,r}^2(s)\gamma_{ii}(s) + D_{j,r}^2(s)\gamma_{jj}(s)\}^{\beta/2}$$

$$\leq Cn^{-\beta/2}m^{\beta/2-1}\sum_{r=1}^m \max_{0\leq s\leq 1}\{E|D_{i,r}^2(s)\gamma_{ii}(s)|^{\beta/2} + |D_{j,r}^2(s)\gamma_{jj}(s)|^{\beta/2}\}$$

$$\leq Cn^{-\beta/2},$$

where the last equality is due to the fact, which was proved in (49)–(50),

$$E\{|D_{j,r}^2(s)\gamma_{jj}(s)|^{\beta/2}\} \leq \sqrt{ED_{j,r}^{2\beta}(s)E\gamma_{jj}^\beta(s)} \leq Cm^{-\beta/2}\max_{0\leq s\leq 1}E\{\gamma_{jj}^\beta(s)\}. \quad \square$$

**Acknowledgments.** The paper was written while the first author was working at NSF. Any conclusions or recommendations expressed in this material are those of the authors and do not necessarily reflect the views of the National Science Foundation. Jian Zou is a Ph. D student at University of Connecticut. The authors thank the editor, associate editor and two anonymous referees for stimulating comments and suggestions, which led to significant improvements of the paper.

DEPARTMENT OF STATISTICS                          DEPARTMENT OF STATISTICS
UNIVERSITY OF WISCONSIN-MADISON                   UNIVERSITY OF CONNECTICUT
MADISON, WISCONSIN 53706                          STORRS, CONNECTICUT 06269
USA                                               USA